\documentclass[twoside,12pt]{article}

\usepackage{amsfonts}
\usepackage{amsmath}

\newcommand{\Xcomment}[1]{}

\newtheorem{theorem}{Theorem}

\newtheorem{corollary}{Corollary}
\newtheorem{proposition}{Proposition}

\newenvironment{myitem}{\refstepcounter{equation}\begin{enumerate}%
\item[(\arabic{equation})]$\quad$}{\end{enumerate}}

\newenvironment{myitem1}{\refstepcounter{equation}\begin{enumerate}%
\item[(\arabic{equation})]\begin{itemize}}%
{\end{itemize}\end{enumerate}}

\makeatletter
\renewcommand{\section}{\@startsection{section}{1}{0pt}%
{-3.5ex plus -1ex minus -.2ex}{2.3ex plus .2ex}%
{\normalfont\Large}}
\makeatother

\def\Rset{{\mathbb R}}
\def\Zset{{\mathbb Z}}
\def\Cset{{\mathbb C}}
\def\Nset{{\mathbb N}}

\def\tilde{\widetilde}

\def\bar{\overline}
\def\eps{\epsilon}

\def\Bscr{{\cal B}}
\def\Cscr{{\cal C}}
\def\Dscr{{\cal D}}
\def\Fscr{{\cal F}}

\def\Hscr{{\cal H}}
\def\Lscr{{\cal L}}
\def\Pscr{{\cal P}}
\def\Rscr{{\cal R}}
\def\Sscr{{\cal S}}

\def\Bstar{{\cal B}^{\ast}}
\def\Estar{E^{\ast}}
\def\Epstar{E^{'\ast}}
\def\Eppstar{E^{''\ast}}
\def\lamblam{\lambda\setminus\bar\lambda}

\def\rest#1{_{\,\vrule height 1.8ex width 0.05em depth 0pt\, #1}}
\newcommand{\refeq}[1]{(\ref{eq:#1})}
\def\qed{ \ \vrule width.2cm height.2cm depth0cm\smallskip}

\def\SC{{\Sscr\Cscr}}
\def\diver{{\rm div}}

\def\trap{\hbox{\unitlength=1mm\begin{picture}(4.5,4)
\put(0,0){{\scriptsize /}}\put(0.5,-0.6){\line(1,0){3.5}}%
\put(1.2,2){\line(1,0){2}}\put(3,0){{\scriptsize $\setminus$}}%
\end{picture}}}

\def\parall{\hbox{\unitlength=1mm\begin{picture}(4.5,4)
\put(0,0){{\scriptsize /}}\put(0.5,-0.6){\line(1,0){3}}%
\put(1.2,2){\line(1,0){3}}\put(3,0){{\scriptsize /}}\end{picture}}}

\begin{document}

\begin{center}
\Large Discrete strip-concave functions,\\
Gelfand-Tsetlin patterns, and related polyhedra
\end{center}

\medskip
\begin{center}
Vladimir~I.~Danilov\footnote[1]
{Central Institute of Economics and Mathematics
of the RAS, 47, Nakhimovskii Prospect, 117418 Moscow, Russia;
emails: vdanilov43@mail.ru, koshevoy@cemi.rssi.ru. Supported in part
by grant NSh-1939.2003.6 from the Russian Foundation of Basic
Research. The third author was also supported by a grant from LIFR
MIIP.
},
Alexander~V.~Karzanov\footnote[2]{Institute for System Analysis of
the RAS, 9, Prospect 60 Let Oktyabrya, 117312 Moscow, Russia;
email: sasha@cs.isa.ru. Corresponding author.
},
and Gleb~A.~Koshevoy$^1$
\end{center}

\bigskip
\begin{quote}
\small {\bf Abstract.} Discrete strip-concave functions considered in
this paper are, in fact, equivalent to an extension of
Gelfand-Tsetlin patterns to the case when the pattern has a
not necessarily triangular but convex configuration. They arise by
releasing one of the three types of rhombus inequalities for discrete
concave functions (or ``hives'') on a ``convex part'' of a triangular
grid. The paper is devoted to a combinatorial study of certain
polyhedra related to such functions or patterns, and results on faces,
integer points and volumes of these polyhedra are presented. Also some
relationships and applications are discussed.

In particular, we characterize, in terms of valid inequalities, the
polyhedral cone formed by the boundary values of discrete
strip-concave functions on a grid having trapezoidal configuration. As
a consequence of this result, necessary and sufficient conditions on a
pair of vectors to be the shape and content of a semi-standard skew
Young tableau are obtained.
 \end{quote}

\noindent
{\em Keywords}: Triangular grid, Gelfand-Tsetlin pattern, Discrete
concave function, Young tableau

\baselineskip 15pt
\parskip 3pt

\section{Introduction} \label{sec:intr}

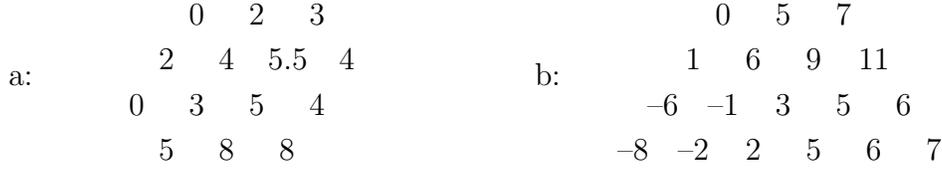
\begin{figure}[htb]
 \begin{center}
  \unitlength=1mm
  \begin{picture}(125,22)
\put(0,10){a:}
\put(20,0){5}
\put(28,0){8}
\put(36,0){8}
\put(16,6){0}
\put(24,6){3}
\put(32,6){5}
\put(40,6){4}
\put(20,12){2}
\put(28,12){4}
\put(34.5,12){5.5}
\put(44,12){4}
\put(24,18){0}
\put(32,18){2}
\put(40,18){3}
\put(70,10){b:}
\put(81,0){--8}
\put(89,0){--2}
\put(98,0){2}
\put(106,0){5}
\put(114,0){6}
\put(122,0){7}
\put(85,6){--6}
\put(93,6){--1}
\put(102,6){3}
\put(110,6){5}
\put(118,6){6}
\put(90,12){1}
\put(98,12){6}
\put(106,12){9}
\put(113,12){11}
\put(94,18){0}
\put(102,18){5}
\put(110,18){7}
  \end{picture}
 \end{center}
 \caption{(a) a hexagonal array with $n=3$, $a=(0,0,0,1)$,
$b=(2,3,3,3)$; (b) a trapezoidal array with $n=3$, $a=(0,0,0,0)$,
$b=(2,3,4,5)$.}
  \label{fig:2arr}
  \end{figure}

Let $n\in\Nset$. Consider a two-dimensional array
$X=(x_{ij})_{0\le i\le n,\,a_i\le j\le b_i}$ of reals, where the index
bounds $a_i,b_i$ (depending on rows) are integers satisfying
$a_i\le b_i$ and:
  \begin{gather}
a_0=0,\quad 0\le a_1-a_0\le a_2-a_1\le \ldots\le a_n-a_{n-1}\le 1,
            \notag\\ \mbox{and}\quad
1\ge b_1-b_0\ge b_2-b_1\ge \ldots\ge b_n-b_{n-1}\ge 0.
                                      \label{eq:conv_arr}
  \end{gather}
We denote the set of pairs $ij$ of indices in $X$ by $V$ and
say that $X$ has {\em convex configuration}.
(This term is justified by the fact that $V$ can be identified with the
set of nodes of a convex triangular grid; see Remark~1 below. We
visualize $X$ so that $(x_{00},\ldots,x_{0b_0})$ is the topmost row
and each triple $x_{ij},x_{i+1,j},x_{i+1,j+1}$ or
$x_{ij},x_{i,j+1},x_{i+1,j+1}$ is disposed so as to form an
equilateral triangle. Then the array is shaped like a convex polygon,
with 3 to 6 sides.) Two examples of such arrays are depicted in
Fig.~\ref{fig:2arr}.

Depending on the shape of the corresponding convex polygon, we may
speak of hexagonal configuration, pentagonal configuration, and etc.
Although main results in this paper will be applicable to any of
these, three special cases with $a_1=\ldots=a_n=0$ are of most
interest for us: (a) $b_i=i$ for each $i$ (giving a
$\Delta$-{\em array}); (b) $b_i=i+m$ for each $i$ (a \trap-{\em array}),
see Fig.~\ref{fig:2arr}b; (c) $b_i=m$ for each $i$
(a \parall-{\em array}), where $m\in\Nset$.
In these cases we will also refer to an array as having {\em
triangular, trapezoidal, {\em or} parallelogram-wise configuration},
respectively (usually ignoring other possible dispositions of triangle,
trapezoid, or parallelogram).
We say that $X$ has {\em
size} $n$ in case (a), and $(n,m)$ in cases (b),(c). Sometimes we will
admit $m=0$ in case (b), regarding $\Delta$-arrays as a degenerate
case of \trap-arrays.

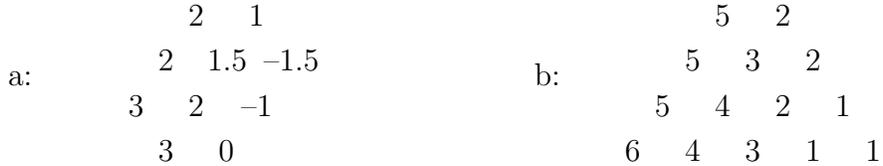
\begin{figure}[htb]
 \begin{center}
  \unitlength=1mm
  \begin{picture}(115,22)
\put(0,10){a:}
\put(20,0){3}
\put(28,0){0}
\put(16,6){3}
\put(24,6){2}
\put(31,6){--1}
\put(20,12){2}
\put(26.5,12){1.5}
\put(34,12){--1.5}
\put(24,18){2}
\put(32,18){1}
\put(70,10){b:}
\put(82,0){6}
\put(90,0){4}
\put(98,0){3}
\put(106,0){1}
\put(114,0){1}
\put(86,6){5}
\put(94,6){4}
\put(102,6){2}
\put(110,6){1}
\put(90,12){5}
\put(98,12){3}
\put(106,12){2}
\put(94,18){5}
\put(102,18){2}
  \end{picture}
 \end{center}
 \caption{Gelfand-Tsetlin pattern examples:
(a) $\partial X$ for $X$ in Fig.~\ref{fig:2arr}a;
(b) $\partial X$ for $X$ in Fig.~\ref{fig:2arr}b.}
  \label{fig:der}
  \end{figure}

Let us associate with $X$ the array $\partial X=
(\partial x_{ij})_{0\le i\le n,\,a_i+1\le j\le b_i}$ of local
differences $\partial x_{ij}:=x_{ij}-x_{i,j-1}$, referring to
$\partial X$ as the {\em row derivative} of $X$. We deal with
arrays $X$ satisfying the following condition:
for $i=1,\ldots,n$ and $j=a_i+1,\ldots,b_i$,
  \begin{equation}  \label{eq:concX}
\partial x_{ij}\ge \partial x_{i-1,j}\;\; \mbox{(when $j\le b_{i-1}$)}
\quad\mbox{and}\quad \partial x_{i-1,j}\ge\partial x_{i,j+1}\;\;
\mbox{(when $j<b_i$)}.
  \end{equation}
The array $\partial X$ obeying~\refeq{concX} and having triangular
configuration is said to be a {\em Gelfand-Tsetlin pattern}, and in
this paper we apply the same name to $\partial X$ with such a
property when $X$ has an arbitrary convex configuration as well.
In this case we call $X$ a {\em strip-concave} array, using
an analogy with the corresponding functions explained in Remark~1
below. For example, both arrays in
Fig.~\ref{fig:2arr} are strip-concave; their row derivatives are shown
in Fig.~\ref{fig:der}.

One can identify the set of all arrays for $V$ with the Euclidean
space $\Rset^V$ whose unit base vectors are indexed by the pairs
$ij\in V$.
Let $\SC_V$ denote the set of arrays $X\in\Rset^V$ that
satisfy property~\refeq{concX} and the normalization condition
$x_{00}=0$; imposing this condition leads to no loss of generality in
what follows. Then $\SC_{V}$ is a polyhedral cone in $\Rset^{V}$.

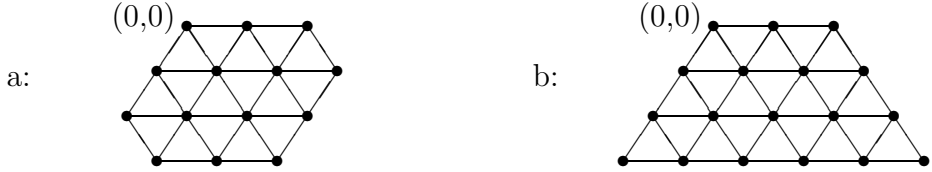
\begin{figure}[htb]
 \begin{center}
  \unitlength=1mm
  \begin{picture}(125,20)
\put(0,10){a:}
\put(20,0){\circle*{1.5}}
\put(28,0){\circle*{1.5}}
\put(36,0){\circle*{1.5}}
\put(16,6){\circle*{1.5}}
\put(24,6){\circle*{1.5}}
\put(32,6){\circle*{1.5}}
\put(40,6){\circle*{1.5}}
\put(20,12){\circle*{1.5}}
\put(28,12){\circle*{1.5}}
\put(36,12){\circle*{1.5}}
\put(44,12){\circle*{1.5}}
\put(24,18){\circle*{1.5}}
\put(32,18){\circle*{1.5}}
\put(40,18){\circle*{1.5}}
\put(20,0){\line(1,0){16}}
\put(16,6){\line(1,0){24}}
\put(20,12){\line(1,0){24}}
\put(24,18){\line(1,0){16}}
\put(16,6){\line(2,3){8}}
\put(20,0){\line(2,3){12}}
\put(28,0){\line(2,3){12}}
\put(36,0){\line(2,3){8}}
\put(20,0){\line(-2,3){4}}
\put(28,0){\line(-2,3){8}}
\put(36,0){\line(-2,3){12}}
\put(40,6){\line(-2,3){8}}
\put(44,12){\line(-2,3){4}}
\put(14,18){(0,0)}
\put(70,10){b:}
\put(82,0){\circle*{1.5}}
\put(90,0){\circle*{1.5}}
\put(98,0){\circle*{1.5}}
\put(106,0){\circle*{1.5}}
\put(114,0){\circle*{1.5}}
\put(122,0){\circle*{1.5}}
\put(86,6){\circle*{1.5}}
\put(94,6){\circle*{1.5}}
\put(102,6){\circle*{1.5}}
\put(110,6){\circle*{1.5}}
\put(118,6){\circle*{1.5}}
\put(90,12){\circle*{1.5}}
\put(98,12){\circle*{1.5}}
\put(106,12){\circle*{1.5}}
\put(114,12){\circle*{1.5}}
\put(94,18){\circle*{1.5}}
\put(102,18){\circle*{1.5}}
\put(110,18){\circle*{1.5}}
\put(82,0){\line(1,0){40}}
\put(86,6){\line(1,0){32}}
\put(90,12){\line(1,0){24}}
\put(94,18){\line(1,0){16}}
\put(82,0){\line(2,3){12}}
\put(90,0){\line(2,3){12}}
\put(98,0){\line(2,3){12}}
\put(106,0){\line(2,3){8}}
\put(114,0){\line(2,3){4}}
\put(90,0){\line(-2,3){4}}
\put(98,0){\line(-2,3){8}}
\put(106,0){\line(-2,3){12}}
\put(114,0){\line(-2,3){12}}
\put(122,0){\line(-2,3){12}}
\put(84,18){(0,0)}
  \end{picture}
 \end{center}
 \caption{(a) the grid for the array in Fig.~\ref{fig:2arr}a;
(b) the grid for the array in Fig.~\ref{fig:2arr}b.}
  \label{fig:grid}
  \end{figure}

\medskip
\noindent
\underline{Remark~1}.
Let $\alpha,\beta$ be linearly independent vectors in $\Rset^2$.
By a {\em convex (triangular) grid} we mean a finite planar graph
$G=(V,E)$ embedded in the plane so that each node of $G$ is a point
with integer coordinates $(i,j)$ in the basis $(\alpha,\beta)$,
each edge is the straight-line segment connecting a pair $u,v$ of
nodes with $u-v\in\{\alpha,\beta,\alpha+\beta\}$, each bounded face
is a triangle with three edges (a {\em little triangle} of $G$),
and the union $\Rscr$ of bounded faces covers all
nodes and forms a convex polygon in the plane. A convex grid can be
considered up to an affine transformation, and to agree with the
above visualization of arrays, one should take the generating vectors
as, e.g., $\alpha:=(-1/2,-\sqrt{3}/2)$ and $\beta:=(1,0)$ and assume
that $(0,0)\in V$ and $(i,j)\ge(0,0)$ for all $(i,j)\in V$.
(The convex grids behind the arrays in Fig.~\ref{fig:2arr} are exposed
in Fig.~\ref{fig:grid}.) A function
$x:V\to\Rset$ determines an array $X$ of convex configuration in a
natural way: $x_{ij}:=x(i,j)$. The arrays in $\SC_V$ (considering
$V$ as the index set) are determined by the functions $x$ having
the following property: if $f$ is the extension of $x$ to $\Rscr$
which is affinely linear on each bounded face of $G$, then $f$ is a
concave function within each region (strip) confined by the boundary
of $G$ and lines $i\alpha+\Rset\beta$ and $(i-1)\alpha+\Rset\beta$,
$i=1,2,\ldots$. We call such a function $x$ {\em discrete
strip-concave} (by an analogy with discrete concave functions; see
Remark~2 in the end of this section), and accordingly apply the
adjective ``strip-concave'' to the arrays with property~\refeq{concX}.

\medskip
Local differences on the ``boundary'' of $X$ will be of most
interest for us in this paper. These are represented by four tuples
$\lambda^X,\bar\lambda^X,\mu^X,\nu^X$ (concerning the lower, upper,
left and right boundaries, respectively) defined by
  \begin{gather}
\lambda^X_j:=\partial x_{nj},\;\; j=1,\ldots,b_n;\qquad
\bar\lambda^X_{j'}:=\partial x_{0j'},\;\; j'=1,\ldots,b_0; \notag\\
\mu^X_i:=x_{ia_i}-x_{i-1,a_{i-1}}\quad \mbox{and}\quad
\nu^X_i:=x_{ib_i}-x_{i-1,b_{i-1}},\quad i=1,\ldots,n. \notag
  \end{gather}
($\bar\lambda^X$ vanishes when $b_0=0$.)
For example, the array $X$ in Fig.~\ref{fig:2arr}a has
$\lambda^X=(3,0)$, $\bar\lambda^X=(2,1)$, $\mu^X=(2,-2,5)$ and
$\nu^X=(1,0,4)$, and
the array $X$ in Fig.~\ref{fig:2arr}b has $\lambda^X=(6,4,3,1,1)$,
$\bar\lambda^X=(5,2)$, $\mu^X=(1,-7,-2)$ and $\nu^X=(4,-5,1)$.

Given $\lambda=(\lambda_{a_n+1},\ldots,\lambda_{b_n})$, $\bar\lambda=
(\bar\lambda_1,\ldots,\bar\lambda_{b_0})$ and $\mu,\nu\in\Rset^n$,
define
  $$
\SC(\lamblam,\mu,\nu):=\{X\in\SC_{V}:
(\lambda^X,\bar\lambda^X,\mu^X,\nu^X)=(\lambda,\bar\lambda,\mu,\nu)\}.
  $$
This set, if nonempty, forms a bounded polyhedron (a polytope)
in $\Rset^{V}$ in case of \trap- and \parall-arrays.
Indeed,~\refeq{concX} and $x_{00}=0$ imply
  \begin{equation} \label{eq:bounds}
x_{ij}\le \mu_1+\ldots+\mu_i+\lambda_1+\ldots+\lambda_j
\qquad\mbox{and}\qquad x_{ij}\ge \mu_1+\ldots+\mu_i+q,
  \end{equation}
where $q:=\lambda_{n-i+1}+\ldots
+\lambda_{n-i+j}$ for \trap-arrays, and $q:=\bar\lambda_1+\ldots+
\bar\lambda_j$ for \parall-arrays. (On the other hand, such a
polyhedron $\Pscr$ is unbounded when there is at least one
interior entry and both left and right boundaries make a bend,
i.e., $0<a_n<n$ and $0<b_n-b_0<n$; in particular, if the hexagonal
configuration takes place. One can check that adding any positive
constant to all interior entries of an array $X\in\Pscr$ gives a point
in $\Pscr$ as well.)

The first problem we deal with in this paper is to characterize
the set $\Bscr_{V}$ of all quadruples $(\lambda,\bar\lambda,\mu,\nu)$
(depending on ${V}$) such that
$\SC(\lamblam,\mu,\nu)$ is nonempty.
Two conditions on such quadruples are trivial. The first one
comes up from the fact that~\refeq{concX} implies that $\lambda^X$ is
weakly decreasing, i.e., $\lambda^X_{a_n+1}\ge\ldots\ge
\lambda^X_{b_n}$, and similarly for $\bar\lambda$.
The second one comes up by observing that
  $$
|\lambda^X|-|\bar\lambda^X|+|\mu^X|-|\nu^X|=(x_{nb_n}-x_{na_n})
-(x_{0b_0}-x_{00})+(x_{na_n}-x_{00})-(x_{nb_n}-x_{0b_0})=0,
  $$
where for a tuple (vector) $d=(d_p,\ldots,d_q)$, $|d|$ stands for
$\sum(d_i: i=p,\ldots,q)$.

To obtain the desired characterization, we need to introduce certain
values depending on $\lambda,\bar\lambda$. For $k\in\Zset_+$,
define
  $$
\delta_k(j):=\max\{0,\bar\lambda_{j-k}-\lambda_j\},\;\;
j=a_n+1,\ldots,b_n, \quad\mbox{and}\quad
\Delta_k:=\delta_k(a_n+1)+\ldots+\delta_k(b_n),
  $$
letting by definition $\delta_k(j):=0$ if $j-k\le 0$ or $j-k>b_0$. We
refer to $\Delta_k$ as the $k$-th {\em deficit} of
$\lamblam$.

We shall explain later that the above problem is reduced to the case of
trapezoidal configuration. Necessary and sufficient conditions on the
corresponding quadruples for \trap-arrays are given in the
following theorem. Hereinafter, for $d=(d_p,\ldots,d_q)$ and $I\subseteq
\{p,\ldots,q\}$, $d(I)$ denotes $\sum(d_i: i\in I)$, and for
$p\le k\le k'\le q$, $d[k,k']$ denotes $d_k+\ldots+d_{k'}$.

  \begin{theorem}  \label{tm:t1}
For $n\in\Nset$ and $m\in\Zset_+$, let $\lambda=(\lambda_1,\ldots,
\lambda_{n+m})$ and $\bar\lambda=(\bar\lambda_1,\ldots,\bar\lambda_m)$
be weakly decreasing, and let $\mu,\nu\in\Rset^n$ be such that
$|\lambda|-|\bar\lambda|+|\mu|-|\nu|=0$. Then a strip-concave
\trap-array $X$ with $(\lambda^X,\bar\lambda^X,\mu^X,\nu^X)=
(\lambda,\bar\lambda,\mu,\nu)$ exists if and only if the
inequality
  \begin{equation}  \label{eq:horn}
\lambda[1,|I|]+\mu(I)-\nu(I)-\Delta_{|I|}\ge 0
  \end{equation}
holds for each (including empty) subset $I\subseteq\{1,\ldots,n\}$.
Furthermore, if $\lambda,\bar\lambda,\mu,\nu$ are integer and the
polytope $\SC(\lamblam,\mu,\nu)$ is nonempty, then
it contains an integer point.
  \end{theorem}

In particular, $\Bscr_{V}$ forms a polyhedral cone
(in $\Rset^{n+m}\times\Rset^m\times\Rset^n\times\Rset^n$) for ${V}$
in question. Also~\refeq{horn} implies evident relations
$\lambda_j\ge \bar\lambda_j$ ($j=1,\ldots,m$) and
$\lambda_j\le\bar\lambda_{j-n}$ ($j=n+1,\ldots,n+m$), where the former
is easily obtained by taking $I=\emptyset$,
and the latter by comparing $|\lambda|-|\bar\lambda|+|\mu|-|\nu|=0$
with~\refeq{horn} for $I=\{1,\ldots,n\}$.

Note that relation~\refeq{horn} involves a piece-wise linear term,
namely, $\Delta_{|I|}$. One can replace each instance of~\refeq{horn}
by a collection of $2^m$ linear inequalities, yielding an equivalent
version of Theorem~\ref{tm:t1}. It turns out that typically these
inequalities determine facets of the cone $\Bscr_V$ for \trap-case;
the precise list of facets of this cone is established in
Section~\ref{sec:facet}.
(We shall see that the number of facets grows exponentially in $n,m$.
On the other hand, to verify that a given quadruple
$(\lambda,\bar\lambda,\mu,\nu)$ belongs to $\Bscr_V$, it suffices to
check validity of~\refeq{horn} only for $n+1$ sets $I$: for
$k=0,\ldots,n$, take $I$ with $|I|=k$ maximizing $(\nu-\mu)(I)$.)

For an arbitrary convex configuration, the problem with prescribed
local differences $\lambda,\bar\lambda,\mu,\nu$ is reduced to the
trapezoidal case as follows. Since the polyhedron
$\Pscr:=\SC(\lamblam,\mu,\nu)$ is described by a linear system formed
by the inequalities in~\refeq{concX} and the corresponding
equalities involving $\lambda,\bar\lambda,\mu,\nu$,
one can efficiently compute a number $c\in\Rset_+$ such that
if $\Pscr$ is nonempty, then there exists $X\in\Pscr$ with
$|x_{ij}|<c/2$ for all entries $x_{ij}$. (For example, one can roughly
take $c$ equal to $|V|^{|V|}$ times the maximum absolute value
$\alpha$ of the entries in $\lambda,\bar\lambda,\mu,\nu$, taking into
account that the constraint matrix of the system has entries 0,1,--1.
In fact, there is a bound $c$ linear in $\alpha|V|$;
cf.~\refeq{bounds} for \parall-arrays.) Suppose $a_n\ne 0$ and take
the maximum $p$ with $a_p=0$ (then $a_i=i-p$ for $p<i\le n$).
Add to ${V}$ the set $A$ of pairs $ij$ with $0\le j<i-p\le n-p$,
define $\lambda'_j:=c$ for $j=1,\ldots,n-p$, and define
$\mu'_i:=\mu_i-c$ for $i=p+1,\ldots,n$. Symmetrically,
if $b_n<b_0+n$, we take the maximum $q$ with $b_q=b_0+q$, add the set
$B$ of pairs $ij$ with $1\le j-b_n\le i-q \le n-q$, define
$\lambda'_j:=-c$ for $j=b_n+1,\ldots,b_n+n-q$, and define
$\nu'_i:=\nu_i-c$ for $i=q+1,\ldots,n$. Let $\lambda'$ coincide with
$\lambda$ for the remaining entries, and similarly for $\mu',\nu'$.
The resulting ${V}':={V}\cup A\cup B$ gives a trapezoid (of size
$(n,b_0)$), and it is straightforward to verify that the set
$\Pscr':=\SC(\lambda'/\bar\lambda,\mu',\nu')$ (concerning ${V}'$) is
nonempty if and only if $\Pscr$ is so, that the restriction of any
$X'\in\Pscr'$ to ${V}$ belongs to $\Pscr$, and that $X$ as above is
extended in a natural way to an array in $\Pscr'$.

Applying this reduction to the parallelogram-wise
configuration of size $(n,m)$, one can derive the following
corollary from Theorem~\ref{tm:t1}.

  \begin{corollary}  \label{cor:par}
Let $n,m\in\Nset$, and let $\mu,\nu\in\Rset^n$ and weakly decreasing
$\lambda,\bar\lambda\in\Rset^m$ satisfy
$|\lambda|-|\bar\lambda|+|\mu|-|\nu|=0$. Then a strip-concave
{\rm \parall}-array $X$ with $(\lambda^X,\bar\lambda^X,\mu^X,\nu^X)=
(\lambda,\bar\lambda,\mu,\nu)$ exists if and only if
for each subset $I\subseteq\{1,\ldots,n\}$, the inequality
  $$
  \lambda[1,|I|]-\bar\lambda[m-|I|+1,m]+\mu(I)-\nu(I)-\Delta_{|I|}\ge 0
  $$
holds for $|I|\le m$, and the inequality
  $$
   |\lambda|-|\bar\lambda|+\mu(I)-\nu(I)\ge 0
  $$
holds for $|I|>m$. Furthermore, if $\lambda,\bar\lambda,\mu,\nu$ are
integer and the polytope $\SC(\lamblam,\mu,\nu)$ is nonempty, then it
contains an integer point.
  \end{corollary}

(To see this, observe that each entry $\nu'_i$ for the new right
boundary tuple is equal to $\nu_i-c$, that $\mu'=\mu$,  and that
$\lambda'_j=-c$ for $j=m+1,\ldots,m+n$. The fact that $\bar\lambda$
has all entries greater than $-c$ implies that for $k=0,\ldots,n$, each
$j$ with $\max\{m,k\}<j\le m+k$ contributes $\bar\lambda_{j-k}+c$ units
to the new $k$-deficit $\Delta'_k$ (whereas $\delta'_k(j)=\delta_k(j)$
for $j=1,\ldots,m$ and $\delta'_k(j)=0$ for the remaining $j$'s).
Therefore, given $I\subseteq\{1,\ldots,n\}$, the new $|I|$-deficit
becomes $\Delta_{|I|}+\bar\lambda[m+1-|I|,m]+|I|c$ whenever
$|I|\le m$, and $|\bar\lambda|+mc$ whenever $|I|>m$. Also
$\lambda'[1,|I|]=\lambda[1,|I|]$ if $|I|\le m$, and
$\lambda'[1,|I|]=|\lambda|-(|I|-m)c$ if $|I|>m$. Now
Corollary~\ref{cor:par} is obtained from Theorem~\ref{tm:t1}
by substituting these relations, together with
$\mu'(I)=\mu(I)$ and $\nu'(I)=\nu(I)-|I|c$, into
relation~\refeq{horn} (taken with primes).)

A converse reduction, from \trap- to \parall-case, is easily
constructed as well, and Theorem~\ref{tm:t1} follows from
Corollary~\ref{cor:par}. In contrast, we cannot point out a ``simple''
reduction of Theorem~\ref{tm:t1} to its special case with $m=0$
concerning $\Delta$-arrays. (Nevertheless, a more intricate, though
constructive, way of reducing does exist, as we explain in part D of
Section~\ref{sec:concl}. In fact, this sort of reduction is behind our
method of proof of Theorem~\ref{tm:t1} where the case $m=0$ is used as
a base.)

\medskip
Another object of our study is the set of vertices of the
polyhedron formed by strip-concave arrays $X$ with convex
configuration whose entries are fixed only on the lower, upper and
left boundaries. More precisely, for $\lambda=(\lambda_{a_n+1},\ldots,
\lambda_{b_n})$, $\bar\lambda=(\bar\lambda_1,\ldots,
\bar\lambda_{b_0})$ and $\mu\in\Rset^n$, define
  $$
\SC(\lamblam,\mu):=\{X\in\SC_{V}:
(\lambda^X,\bar\lambda^X,\mu^X)=(\lambda,\bar\lambda,\mu)\}.
  $$
(This polyhedron is bounded in case of $\Delta$-, \trap-, or
\parall-configuration since the bounds on $x_{ij}$ indicated
in~\refeq{bounds} remain valid in this case too.) We show the
following.

   \begin{theorem}  \label{tm:t2}
For an arbitrary convex configuration and integer
$\lambda,\bar\lambda,\mu$, the polyhedron
$\SC(\lamblam,\mu)$ is integral, i.e., each face of
this polyhedron contains an integer point.
   \end{theorem}

Note that for arbitrary reals $q_1,\ldots,q_n$,
the transformation of an array $X$ into the array $X'$ with
entries $x'_{ij}:=x_{ij}+q_i$ preserves the row derivative.
Such a transformation
shifts a polyhedron $\SC(\lamblam,\mu,\nu)$ into
$\SC(\lamblam,\mu',\nu')$ with
$\mu'_i:=\mu_i+q_i-q_{i-1}$ and $\nu'_i:=\nu_i+q_i-q_{i-1}$ (letting
$q_0:=0$) and it maintains relation~\refeq{horn}. This implies
that, without loss of generality, in Theorem~\ref{tm:t1} one can
consider only the quadruples of the form
$(\lambda,\bar\lambda,0^n,\nu)$ (where $0^n$ is the zero $n$-tuple).
Similarly, one can restrict $\mu$ to be $0^n$ in Theorem~\ref{tm:t2}
as well.

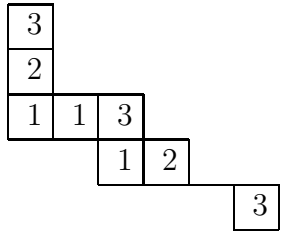
\begin{figure}[htb]
 \begin{center}
  \unitlength=1mm
  \begin{picture}(40,30)
\put(30,0){\line(1,0){6}}
\put(12,6){\line(1,0){24}}
\put(0,12){\line(1,0){24}}
\put(0,18){\line(1,0){18}}
\put(0,24){\line(1,0){6}}
\put(0,30){\line(1,0){6}}
\put(0,12){\line(0,1){18}}
\put(6,12){\line(0,1){18}}
\put(12,6){\line(0,1){12}}
\put(18,6){\line(0,1){12}}
\put(24,6){\line(0,1){6}}
\put(30,0){\line(0,1){6}}
\put(36,0){\line(0,1){6}}
\put(32.5,2){3}
\put(14.5,8){1}
\put(20.5,8){2}
\put(2.5,14){1}
\put(8.5,14){1}
\put(14.5,14){3}
\put(2.5,20){2}
\put(2.5,26){3}
  \end{picture}
 \end{center}
 \caption{the semi-standard skew Young tableau
corresponding to the pattern in Fig.~\ref{fig:der}b (here $\lambda=
(6,4,3,1,1)$, $\bar\lambda=(5,2)$ and $\nu=(3,2,3)$).}
  \label{fig:tabl}
  \end{figure}

When dealing with $\Delta$-configuration, for a triple
$(\lambda,0^n,\nu)$, inequality~\refeq{horn} turns into
the {\em majorization condition} $\lambda[1,|I|]\ge \nu(I)$. Therefore,
for a fixed $\lambda$, the set $\{\nu: (\lambda,0^n,\nu)\in\Bscr_{V}\}$
forms a {\em permutohedron}, a polytope
$\Pscr$ formed by all vectors $z\in\Rset^n$ with the same value $|z|$
such that for $k=1,\ldots,n-1$, the sum of any $k$ entries of $z$ does
not exceed a constant depending only on $k$. (The vertices of $\Pscr$
are obtained by permuting entries of a fixed $n$-vector $h$; in our
case, $h=\lambda$.) It is known that for
nonnegative integer $\lambda,\nu$, the majorization
condition is necessary and sufficient for the existence of a
{\em semi-standard Young tableau with shape $\lambda$
and content $\nu$}, and that these tableaux one-to-one correspond to
the integer Gelfand-Tsetlin patterns respecting $\lambda,\nu$; for a
definition and a survey, see~\cite{St}.
Theorem~\ref{tm:t1} (Corollary~\ref{cor:par}) shows that in case of
\trap-arrays (resp. \parall-arrays) and $\lambda,\bar\lambda$ fixed,
the analogous set $\{\nu: (\lambda,\bar\lambda,0^n,\nu)\in
\Bscr_{V}\}$ forms a permutohedron in $\Rset^n$ as well (but now the
corresponding vertex generating vector $h$ becomes less trivial to
write down; it will be indicated in Section~\ref{sec:concl}). Each
integer (generalized) Gelfand-Tsetlin pattern for nonnegative integer
$\lambda,\bar\lambda,\nu$ determines a so-called
semi-standard {\em skew} Young tableau {\em with shape $\lamblam$ and
content} $\nu$ (cf.~\cite{St}), and our theorem (corollary) yields
necessary and sufficient conditions for the existence of such
tableaux. Figure~\ref{fig:tabl} illustrates an instance of
semi-standard skew Young tableau.

It should be noted that in case of $\Delta$-configuration one can
obtain the claim of Theorem~\ref{tm:t2} by using a description for the
generators of the Gelfand-Tsetlin patterns cone given in Berenstein
and Kirillov~\cite{BK}.

Our method of proof of Theorem~\ref{tm:t2} is based on attracting
a certain equivalent flow model and showing that the integer points in
$\SC(\lamblam,0^n)$ one-to-one correspond to the
integer flows in a certain directed graph. In addition, we explain how
to use the flow approach to easily show that Kostka coefficient
$K(\lambda,\nu)$ (or $K(\lamblam,\nu)$), as well as the intrinsic
volume of $\SC(\lambda,0^n,\nu)$ (resp. $\SC(\lamblam,0^n,\nu)$) in
the nondegenerate case, preserves under a permutation of the entries
of $\nu$. Here $K(\lambda,\nu)$ is the number of semi-standard Young
tableaux with shape $\lambda$ and content $\nu$ (which is equal to the
number of integer points in $\SC(\lambda,0^n,\nu)$), while
$K(\lamblam,\nu)$ concerns the corresponding skew tableaux.

This paper is organized as follows.
Theorems~\ref{tm:t1} and~\ref{tm:t2} are proved
in Sections~\ref{sec:proof1} and~\ref{sec:proof2}, respectively.
Section~\ref{sec:facet} is devoted to a sharper version of
Theorem~\ref{tm:t1} that precisely describes the set of linear
inequalities determining facets of the cone $\Bscr_V$ in \trap-case
(Theorem~\ref{tm:t3}). The concluding Section~\ref{sec:concl}
discusses some additional aspects related to these theorems and
demonstrates consequences from the proving method of
Theorem~\ref{tm:t2}: a combinatorial characterization of the vertices
of polyhedra $\SC(\lamblam,\mu)$, the above-mentioned facts on integer
points and volumes, and others.

We conclude this section with two more remarks.

\medskip
\noindent
\underline{Remark~2}.
Let us say that an array $X$ (as in~\refeq{conv_arr}) is {\em (fully)
concave} if it satisfies~\refeq{concX} and
  \begin{equation}  \label{eq:add}
  x_{ij}-x_{i+1,j}\ge x_{i-1,j-1}-x_{i,j-1} \qquad \mbox{for all
$1\le i<n$ and $a_i<j\le b_i$}.
  \end{equation}
This is equivalent to saying that the extension $f$ of the
function $x$ on the nodes of the corresponding grid $G$ (cf.
Remark~1) is concave in the entire region $\Rscr$. The functions $x$
with such a property are often called {\em discrete concave} ones,
and a series of interesting results on these have been obtained.
Knutson, Tao and Woodward~\cite{KTW} pointed out the precise list of
facets of the cone ${\rm BNDR}_n$ formed by all possible triples
$(\lambda,\mu,\nu)$ of $n$-tuples whose entries are the differences
$x(v)-x(u)$ on boundary edges $uv$ for a discrete concave function $x$
on the triangular grid of size $n$, or a {\em hive} (equivalently:
$\lambda,\mu,\nu$ are the spectra of three Hermitian $n\times n$
matrices with zero sum). Also it is shown
in~\cite{KT} that for each integer $(\lambda,\mu,\nu)\in {\rm BNDR}_n$
there exists an integer discrete concave function $x$ as required for
this triple. (A history of studying this cone and related topics are
reviewed in~\cite{Ful}, see also~\cite{DK}). Nontrivial constraints
for ${\rm BNDR}_n$ are expressed by Horn's inequalites. These are
generalized to an arbitrary convex grid (see~\cite{Kar-03}), and
relation~\refeq{horn} in Theorem~\ref{tm:t1} is, in essense,
equivalent to a special case of Horn's inequalites. We will briefly
explain in Section~\ref{sec:concl} that Theorem~\ref{tm:t1} can be
derived from the above-mentioned results on discrete concave
functions. At the same time, our direct proof of Theorem~\ref{tm:t1}
is much simpler compared with the proofs of the corresponding theorems
in~\cite{KT,KTW}.

\medskip
\noindent
\underline{Remark~3}.
The polyhedron integrality claimed in
Theorem~\ref{tm:t2} need not hold when the array entries are fixed on
the whole boundary. More precisely, by a result due to De~Loera and
McAllister~\cite{DM}, for any $k\in\Nset$, there exist
$\lambda,\mu,\nu\in\Zset^n$ and a triangular array $X$ of size $n$,
with $n=O(k)$, such that $X$ is a vertex of the polytope
$\SC(\lambda,\mu,\nu)$ and some entry of $X$ has denominator $k$.
(Some ingredient from a construction in~\cite{DM} is borrowed
by~\cite{Kar-04} to obtain an analogous result for fully concave
triangular arrays in the case when the values are fixed
only on two ``sides''.) Nevertheless, for \trap-, \parall- or
$\Delta$-configuration, at least one integer vertex in each nonempty
polytope $\SC(\lamblam,\mu,\nu)$ with
$\lambda,\bar\lambda,\mu,\nu$ integer does exist, as explained in
the end of Section~\ref{sec:concl}.

\section{Proof of Theorem~\ref{tm:t1}} \label{sec:proof1}

As explained in the Introduction, it suffices to consider the case
$\mu=0^n$.

To show part ``only if'' in the theorem, we use induction on $n$. Case
$n=1$ is trivial, so assume $n>1$. Let $(\lambda,\bar\lambda,0^n,\nu)
\in\Bscr_{V}$ (for ${V}$ determined by $n,m$) and consider an array
$X\in\SC(\lamblam,0^n,\nu)$ and a set
$I=\{i(1),\ldots,i(k)\}$ with $1\le i(1)<\ldots< i(k)\le n$.

Define $I':=I\cap\{1,\ldots,n-1\}$ and $\lambda'_j:=
\partial x_{n-1,j}$ for $j=1,\ldots,n+m-1$. Then $\lambda_j\ge
\lambda'_j\ge\lambda_{j+1}$ (by~\refeq{concX}). By induction,
  \begin{equation}  \label{eq:phorn}
\lambda'[1,|I'|]-\nu(I')-\Delta'_{|I'|}\ge 0,
  \end{equation}
where $\Delta'_{k'}$ stands for the $k'$-th deficit for
$\lambda',\bar\lambda$, i.e., $\Delta'_{k'}:=\delta'_{k'}(1)+
\ldots+\delta'_{k'}(n+m-1)$, where $\delta'_{k'}(j):=\max\{0,
\bar\lambda_{j-k'}-\lambda'_j\}$. Two cases are possible.

\smallskip
{\em Case 1.} Let $n\not\in I$, i.e., $I'=I$. Since
$\delta_k(j)=\max\{0,\bar\lambda_{j-k}-\lambda_j\}$,
$\delta'_k(j)=\max\{0,\bar\lambda_{j-k}-\lambda'_j\}$ and $\lambda_j
\ge \lambda'_j$, we have $\delta_k(j)\le\delta'_k(j)$, implying
$\Delta_k\le\Delta'_k$. Now, using~\refeq{phorn},
   $$
\lambda[1,k]-\nu(I)-\Delta_k\ge \lambda'[1,k]-\nu(I)-\Delta'_k\ge 0,
   $$
and~\refeq{horn} follows (with $\mu=0^n$).

\smallskip
{\em Case 2.} Let $n\in I$. Then $|I'|=k-1$. Summing up~\refeq{phorn}
and the evident equality $|\lambda|-|\lambda'|-\nu_n=0$, we obtain
  \begin{equation}  \label{eq:pphorn}
 \lambda[1,k]+\sum\left( \lambda_j-\lambda'_{j-1}: j=k+1,\ldots,n+m
\right)-\nu(I)-\Delta'_{k-1}\ge 0.
  \end{equation}

Note also that
$\lambda_j+\delta_k(j)=\max\{\lambda_j,\bar\lambda_{j-k}\}$
(in view of $\delta_k(j)=\max\{0,\bar\lambda_{j-k}-\lambda_j\}$),
and similarly $\lambda'_{j-1}+\delta'_{k-1}(j-1)=
\max\{\lambda'_{j-1},\bar\lambda_{j-k}\}$. Since
$\lambda_j\le\lambda'_{j-1}$, we have
$\lambda_j+\delta_k(j)\le\lambda'_{j-1}+\delta'_{k-1}(j-1)$.
Therefore, $\sum (\lambda_j-\lambda'_{j-1}: j=k+1,\ldots,n+m)
-\Delta'_{k-1}\le -\Delta_k$. This together with~\refeq{pphorn}
implies~\refeq{horn}.

\medskip
Next we show part ``if'' in the theorem. We first consider case $m=0$
(i.e., $\Delta$-configuration); in this case all deficits $\Delta_k$
are zeros, which simplifies the consideration.
We use induction on $n$; case $n=1$ is trivial. Let $n>1$ and
let~\refeq{horn} hold for all $I$. In particular,
$\lambda_1-\nu_n\ge 0$ (by taking $I:=\{n\}$). Also, subtracting
inequality~\refeq{horn} with $I=\{1,\ldots,n-1\}$ from the equality
$|\lambda|-|\nu|=0$, we obtain $\lambda_n-\nu_n\le 0$. Therefore,
as $\lambda$ is weakly decreasing, there exists $p\in\{1,\ldots,n-1\}$
such that
  \begin{equation}  \label{eq:lamp}
  \lambda_p\ge\nu_n\quad \mbox{and}\quad
  \lambda_{p+1}\le\nu_n.
  \end{equation}

Assign the $(n-1)$-tuple $\lambda'$ by the following rule:
  \begin{gather}
 \lambda'_j:=\lambda_j\quad\mbox{for $j=1,\ldots,p-1$}; \qquad
 \lambda'_j:=\lambda_{j+1}\quad\mbox{for $j=p+1,\ldots,n-1$}; \notag\\
 \mbox{and}\quad\lambda'_p:=\lambda_p+\lambda_{p+1}-\nu_n.
                                        \label{eq:lamj}
  \end{gather}

Consider the triple $(\lambda',0^{n-1},\nu')$, where
$\nu':=(\nu_1,\ldots,\nu_{n-1})$.
We assert that
  \begin{equation}  \label{eq:lam1I}
 \lambda'[1,|I'|]\ge\nu(I')
  \end{equation}
holds for each $I'\subseteq\{1,\ldots,n-1\}$. Consider two cases,
letting $k:=|I'|$.

\smallskip
(i) Let $k<p$. Then $\lambda'[1,k]=\lambda[1,k]$, and~\refeq{lam1I}
follows from~\refeq{horn} for $I:=I'$.

\smallskip
(ii) Let $k\ge p$. Define $I:=I'\cup\{n\}$. Then
$\lambda'[1,k]=\lambda[1,k+1]-\nu_n$ (by~\refeq{lamj}), and
we have (using~\refeq{horn})
   $$
 \lambda'[1,k]-\nu(I')=\lambda[1,k+1]-\nu_n-\nu(I')
    =\lambda[1,|I|]-\nu(I)\ge 0.
   $$

Thus, \refeq{lam1I} holds for each $I'$.
Also~\refeq{lamp} and~\refeq{lamj} imply $\lambda_j\ge\lambda'_j\ge
\lambda_{j+1}$ for $j=1,\ldots,n-1$ (in particular, $\lambda'$ is
weakly decreasing), and~\refeq{lamj} together with $|\lambda|=|\nu|$
implies $|\lambda'|=|\nu'|$. By induction there exists a strip-concave
$\Delta$-array $X'$ of size $n-1$ with
$(\lambda^{X'},\mu^{X'},\nu^{X'})=(\lambda',0^{n-1},\nu')$.
Assign $x_{ij}:=x'_{ij}$ for $0\le j\le i\le n-1$ and
$x_{nj}:=\lambda[1,j]$ for $1\le j\le n$. The resulting array $X$ of
size $n$ satisfies~\refeq{concX} and has the desired local differences
on the ``sides'', namely, $(\lambda^X,\mu^X,\nu^X)=(\lambda,0^n,\nu)$.
Hence $(\lambda,0^n,\nu)\in\Bscr_{V}$. Also when $\lambda,\nu$ are
integer, the tuple $\lambda'$ defined by~\refeq{lamj} is integer as
well, and the last claim in the theorem (for $m=0$) follows by
induction, as the integrality of $X'$ implies that for $X$.

\medskip
It remains to prove part ``if'' when $m>0$.
Notice that the triple $\lambda,\bar\lambda,\nu$ can be considered up
to adding a constant to all entries (which matches adding a constant
to the array row derivative), so one may assume that $\lambda$ is
nonnegative. Also, by compactness and scaling, w.l.o.g. one may assume
that $\lambda,\bar\lambda,\nu$ are integer (this slightly simplifies
technical details).

We proceed by induction on $m+|\lambda|$; case $|\lambda|=0$ is
trivial. Let~\refeq{horn} hold for all $I$. In particular,
$\lambda_j\ge\bar\lambda_j\ge\lambda_{j+n}$ for $j=1,\ldots,m$. If
$\lambda_{n+m}=\bar\lambda_m$, we make a simple reduction to
\trap-configuration of size $(n,m-1)$ by truncating the tuples
$\lambda,\bar\lambda$ to $\lambda':=(\lambda_1,\ldots,
\lambda_{n+m-1})$ and $\bar\lambda':=(\bar\lambda_1,\ldots,
\bar\lambda_{m-1})$, respectively. (This maintains~\refeq{horn}, and
if $X$ is a required array of size $(n,m-1)$ for
$\lambda',\bar\lambda',\nu$, then adding to $X$ the elements
$x_{i,i+m}:=x_{i,i+m-1}+\bar\lambda_m$ for $i=0,\ldots,n$ produces a
required array of size $(n,m)$ for $\lambda,\bar\lambda,\nu$). A
similar reduction (discarding $\lambda_1,\bar\lambda_1$) is applied
when $\lambda_1=\bar\lambda_1$.

Therefore, one may assume $\lambda_1>\bar\lambda_1$ and
$\bar\lambda_m>\lambda_{n+m}$. Then there are $1\le r\le n+m-1$ and
$1\le s\le m$ such that
  \begin{equation} \label{eq:lamrs}
\lambda_r\ge\bar\lambda_1=\ldots=\bar\lambda_s>\bar\lambda_{s+1}
\qquad\mbox{and}\qquad \bar\lambda_s>\lambda_{r+1},
  \end{equation}
letting $\bar\lambda_{m+1}:=0$. Note that $\lambda_r>
\bar\lambda_{s+1}$ implies $r\le s+n$ and $\bar\lambda_s>\lambda_{r+1}$
implies $r\ge s$. Define
\begin{equation}\label{eq:lpj}
   \lambda'_j := \left\{
    \begin{array}{rl}
  \lambda_j-1, & \qquad j=r-s+1,\ldots,r,\\
  \lambda_j,  & \qquad j=1,\ldots,r-s,r+1,\ldots,n+m;
  \end{array}
        \right.
   \end{equation}
\begin{equation}\label{eq:blpj}
   \bar\lambda'_j := \left\{
    \begin{array}{rl}
  \bar\lambda_j-1, & \qquad j=1,\ldots,s,\\
  \bar\lambda_j,  & \qquad j=s+1,\ldots,m.
  \end{array}
        \right.
   \end{equation}

Then $\lambda',\bar\lambda'$ are weakly decreasing and $|\lambda'|-
|\bar\lambda'|-|\nu|=0$.
We assert that for any $I\subseteq\{1,\ldots,n\}$ and $k:=|I|$:
  \begin{equation}  \label{eq:lnuD}
\lambda'[1,k]-\nu(I)-\Delta'_k\ge 0,
  \end{equation}
denoting by $\Delta'_k$ the $k$-deficit for $\lambda',\bar\lambda'$,
i.e., the sum of numbers $\delta'_k(j):=\max\{0,\bar\lambda'_{j-k}
-\lambda'_j\}$ over $j$. To see this, first of all observe that
$\delta'_k(j)=\delta_k(j)=0$ if $1\le j\le r$ (since $\lambda_j\ge
\bar\lambda_1$ and $\lambda'_j\ge\bar\lambda'_1$,
by~\refeq{lamrs},\refeq{lpj},\refeq{blpj}). Consider three cases.

\smallskip
(a) Let $k\le r-s$. Then for $j=r+1,\ldots,n+m$, we have
$\lambda'_j=\lambda_j$ and $\lambda'_{j-k}=\bar\lambda_{j-k}$ (in view
of $j-k>s$). Hence $\Delta'_k=\Delta_k$. Also $\lambda'[1,k]=
\lambda[1,k]$. Then~\refeq{lnuD} follows from~\refeq{horn}.

\smallskip
(b) Let $r-s<k\le r$. Then $\delta'_k(j)=\delta_k(j)-1$ for
$j=r+1,\ldots,k+s$ (as $1\le j-k\le s$ implies
$\bar\lambda'_{j-k}=\bar\lambda_{j-k}-1\ge \lambda_j=\lambda'_j$), and
$\delta'_k(j)=\delta_k(j)$ for $j=k+s+1,\ldots,n+m$. So $\Delta'_k=
\Delta_k-(k+s-r)$. Also $\lambda'[1,k]=\lambda[1,k]-(k+s-r)$,
and~\refeq{lnuD} follows.

\smallskip
(c) Let $r<k\le n$. Then $\delta'_k(j)=\delta_k(j)-1$ for $j=k+1,
\ldots,k+s$, and $\delta'_k(j)=\delta_k(j)$ for $j=k+s+1,\ldots,n+m$.
So $\Delta'_k=\Delta_k-s$. Also $\lambda'[1,k]=\lambda[1,k]-s$,
and~\refeq{lnuD} follows.

\smallskip
Since $|\lambda'|<|\lambda|$, by induction the set
$\SC(\lambda'\setminus\bar\lambda',0^n,\nu)$ is nonempty and contains
an integer member $X'$. We transform $X'$ into the desired array $X$
for $\lambda,\bar\lambda,\nu$ as follows. Let $\alpha:=
\lambda'_{r-s+1}$ (=$\lambda_{r-s+1}-1$). For $i=0,\ldots,n$, define
$p(i)$ to be the maximum $j$ such that $\partial x'_{ij}>\alpha$,
letting by definition $\partial x'_{i0}:=\infty$. Then
$p(0)=0$, $p(n)=r-s$ and $p(i)\le i$ for each $i$ (as
$\lambda'_{r-s}>\lambda'_{r-s+1}\ge\bar\lambda'_1\ge
\partial x'_{i,i+1}$). For $i=0,\ldots,n$, define
\begin{equation}\label{eq:xij}
   x_{ij}:= \left\{
    \begin{array}{rl}
  x'_{ij}, & \qquad j=0,\ldots,p(i),\\
  x'_{ij}+j-p(i), & \qquad j=p(i)+1,\ldots,p(i)+s, \\
  x'_{ij}+s,  & \qquad j=p(i)+s+1,\ldots,i+m.
  \end{array}
        \right.
   \end{equation}

Observe that $\lambda^X=\lambda$, $\bar\lambda^X=\bar\lambda$ and
$\nu^X=\nu$ (since $x_{i,i+m}=x'_{i,i+m}+s$ for each $i$). Also $X$
satisfies~\refeq{concX}. To see the latter, let $\eps_{ij}:=
\partial x_{ij}-\partial x'_{ij}$ for all corresponding $i,j$;
then $\eps_{ij}\in\{0,1\}$. Using the definition of
$\alpha,p(0),\ldots,p(n)$, relation~\refeq{xij} and the fact that
$X'$ is strip-concave, it is not difficult to conclude that
$\eps_{ij}<\eps_{i-1,j}$ is possible only if $j=p(i)=p(i-1)+1$.
In this case we have
$\partial x'_{ij}\ge \alpha+1>\partial x'_{i-1,j}$,
whence $\partial x_{ij}\ge \partial x_{i-1,j}$.
Similarly, one can see that if $\eps_{i-1,j}<\eps_{i,j+1}$, then
$j=p(i-1)=p(i)$; in this case $\partial x_{i-1,j}\ge
\partial x_{i,j+1}$ follows from $\partial x'_{i-1,j}\ge \alpha+1>
\partial x'_{i,j+1}$. This implies that $X$ is strip-concave.

This completes the proof of Theorem~\ref{tm:t1}.

\medskip
\noindent
\underline{Remark 4}. The $\Delta$-array $X$ recursively constructed
in the second part of the proof is, in fact, a vertex of
the polytope $\SC(\lambda,0^n,\nu)$. This can be seen as follows.
Given $X'\in\SC(\lambda,0^n,\nu)$, let $Q(X')$ be the set of all
equalities of the form $\partial x'_{ij}=\partial x'_{i-1,j}$ or
$\partial x'_{i-1,j}=\partial x'_{i,j+1}$. A trivial observation is
that $X'$ is a vertex of $\SC(\lambda,0^n,\nu)$ if and only if $X'$
is determined by $Q(X')$, i.e., there is no other point $X''$ in this
polytope such that $Q(X'')\supseteq Q(X')$. In our case, the
equalities as in~\refeq{lamj} (in the recursive process)
give the corresponding equalities for $\partial X$; clearly
the latter equalities determine $X$ uniquely, so $X$ is a vertex
of $\SC(\lambda,0^n,\nu)$. Moreover, if $\lambda,\nu$ are integer,
then $X$ is integer as well. This strengthens the last claim
in the theorem for case $m=0$. On the other hand, the construction of
\trap-array $X$ in the third part of the proof does not guarantee that
this $X$ is a vertex of $\SC(\lamblam,0^n,\nu)$. (Although an integer
vertex in this polytope with $\lambda,\bar\lambda,\nu$ integer does
exist, as explained in Section~\ref{sec:concl}.)

\medskip
\noindent
\underline{Remark 5}. One can accelerate the process of constructing
a required \trap-array $X$ in the third part of the proof. Given
(not necessary integer) $\lambda,\bar\lambda,\nu$, define $\rho:=
\bar\lambda_1-\max\{\lambda_{r+1},\bar\lambda_{s+1}\}$, for $r,s$ as
in~\refeq{lamrs}. When $\lambda_1>\bar\lambda_1$ and $\bar\lambda_m>
\lambda_{n+m}$, we can reduce the corresponding entries of $\lambda,
\bar\lambda$ just by $\rho$ (rather than by one), by setting
$\lambda'_j:=\lambda_j-\rho$ and $\bar\lambda'_j:=\bar\lambda_j-\rho$
in the first lines of~\refeq{lpj} and~\refeq{blpj}, respectively (one
shows that~\refeq{horn} is maintained). Given an array $X'$ for
$\lambda',\bar\lambda',\nu$, we iteratively transform $X'$ into an
array for $\lambda,\bar\lambda,\nu$. More precisely, at the first
iteration, for $\alpha,p(0),\ldots,p(n)$ defined as above, we increase
the entries $x'_{ij}$ for $ij$ as in the second and third lines
of~\refeq{xij} by $\eps(j-p(i))$ and by $\eps s$, respectively,
where $\eps$ is the maximum value not exceeding $\rho$ and such that
the resulting array is still strip-concave ($\eps$ is computed
efficiently). If $\eps<\rho$, we apply a similar procedure (at
the second iteration) to the updated $X'$ and $\rho:=\rho-\eps$, and
so on. One shows that after $O(n^2)$ iterations we get $\rho=0$, and
that the final $X'$ is the desired array $X$ for $\lambda,\bar\lambda,
\nu$. Hence the number of operations in the whole process of finding a
member of $\SC(\lamblam,0^n,\nu)$ is polynomial in
$n$. Such a transformation $X'\to X$ is closely related to a
rearrangement of flows (associated with strip-concave arrays) explained
in part D of Section~\ref{sec:concl}.

\section{Facets of the Cone $\Bscr_V$}
\label{sec:facet}

As mentioned in the Introduction, Theorem~\ref{tm:t1} admits a
reformulation in which the piece-wise linear constraints are replaced
by linear ones. More precisely, one can see that for each
$I\subseteq\{1,\ldots,n\}$, inequality~\refeq{horn} is equivalent
to the set of linear inequalities
  \begin{equation}  \label{eq:modhorn}
\lambda[1,|I|]+\lambda(J+|I|)-\bar\lambda(J)+\mu(I)-\nu(I)\ge 0,
  \end{equation}
where $J$ ranges all subsets of $\{1,\ldots,m\}$, and for $k\in\Zset$,
$J+k$ stands for the set $\{j+k: j\in J\}$. In turns out that, as a
rule, each of the latter inequalities is essential, i.e., determines a
facet of the cone $\Bscr_V$. Note also that the description of this
cone given in Theorem~\ref{tm:t1} involves the ``chamber inequalities''
$\lambda_ j\ge\lambda_{j+1}$ and $\bar\lambda_j\ge\bar\lambda_{j+1}$,
so a priori such inequalities may determine facets as well.
The precise list of facets of $\Bscr_V$ is indicated in the following
theorem.

  \begin{theorem} \label{tm:t3}
For \trap-configuration of size $(n,m)$ with $n\ge 1$ and $m\ge 0$,
inequality~\refeq{modhorn} determines a facet of $\Bscr_V$ if and only
if $0<|I|+|J|<n+m$ and either (i)
$|I|\ne 0,n$ (and $J$ is arbitrary), or (ii) $|I|=0$ and $|J|=1$, or
(iii) $|I|=n$ and $|J|=m-1$. Furthermore, all these facets are
different and $\Bscr_V$ has no other facets if $n=1$ or if $n=2$ and
$m=0$. Otherwise the remaining facets are exactly those determined by
the chamber inequalities $\lambda_ j\ge\lambda_{j+1}$ for
$j=1,\ldots,n+m-1$, and $\bar\lambda_{j'}\ge\bar\lambda_{j'+1}$ for
$j'=1,\ldots,m-1$.
In particular, $\Bscr_V$ has $(2^n-2)2^m+n+4m-2$ facets in case
$n\ge 2$, and $2m$ facets in case $n=1$.
  \end{theorem}

\noindent {\bf Proof.}
It is convenient to consider the reduced cone $\Bstar=\Bstar_{n,m}$
formed by the vectors
$(\lambda,\bar\lambda,\nu)\in\Rset^{n+m}\times\Rset^m\times\Rset^n$
such that $(\lambda,\bar\lambda,0^n,\nu)$ belongs to the cone
$\Bscr:=\Bscr_V$.
Since $(\lambda,\bar\lambda,\mu,\nu)\in \Bscr$ if and only if
$(\lambda,\bar\lambda,0^n,\nu-\mu)\in \Bscr$ (as explained in the
Introduction), the cone $\Bstar$ lies in the hyperplane $\Hscr$
defined by $|\lambda|-|\bar\lambda|-|\nu|=0$ and is described by the
inequalities
  \begin{equation}  \label{eq:redhorn}
\lambda[1,|I|]+\lambda(J+|I|)-\bar\lambda(J)-\nu(I)\ge 0
  \end{equation}
for all $I\subseteq\{1,\ldots,n\}$ and $J\subseteq\{1,\ldots,m\}$,
and by the chamber inequalities. Moreover, there is a natural
bijection between the facets of $\Bscr$ and $\Bstar$, namely: for
any $I,J$, inequality~\refeq{modhorn} determines a facet of $\Bscr$ if
and only if~\refeq{redhorn} determines a facet of $\Bstar$, and
similarly for the chamber inequalities.

In light of these observations, our goal is to characterize those of
the instances of~\refeq{redhorn} and of the chamber inequalities that
determine the facets of $\Bstar$. We will argue in terms of the grid
$G=(V,E)$ (defined in Remark 1 in the Introduction). Let
$\Estar=\Estar_{n,m}$ be the set of ``horizontal'' edges $e_{ij}:=
\{(i,j-1),(i,j)\}$ of $G$. An edge $e_{ij}$ with $i=n$ ($i=0$)
is also denoted by $e_j$ (resp. $\bar e_j$); we refer to the sets
$L:=\{e_1,\ldots,e_{n+m}\}$ and $\bar L:=\{\bar e_1,\ldots,
\bar e_{n+m}\}$ as the {\em lower} and {\em upper} boundaries,
respectively. Besides, we will deal with the edges $r_1,\ldots,r_n$
on the {\em right} boundary $R$, where $r_i$ connects the vertices
$(i-1,i+m-1)$ and $(i,i+m)$. We associate with an edge $e\in \Estar$
the function $\chi^e$ on $\Estar$ taking value 1 on $e$, and 0 on the
remaining edges.

Let $h:\Estar\to\Rset$. The {\em border} of $h$ is defined to be the
function $\sigma=\sigma^h$ on $L\cup\bar L\cup R$
coinciding with $h$ on
$L\cup\bar L$ and taking the value $(h(e_{i,1})+\ldots+h(e_{i,i+m}))-
(h(e_{i-1,1})+\ldots+h(e_{i-1,i+m-1}))$ on $r_i$, $i=1,\ldots,n$; we
also use vector notation, identifying $\sigma$ with the corresponding
triple $(\lambda^h,\bar\lambda^h,\nu^h)$,
where $\lambda^h_j:=\sigma(e_j)$,
$\bar\lambda^h_j:=\sigma(\bar e_j)$ and $\nu^h_i:=\sigma(r_i)$.
Also $h$ is identified with the corresponding array $\partial X$ of
row derivatives (i.e., $\partial x_{ij}=h(e_{ij})$ for all $ij$), and
we say that $h$ is a Gelfand-Tsetlin pattern, or, briefly, a
{\em GT-pattern}, if $\partial X$ is such. In other words, $h$ is a
GT-pattern if $h(e)\ge h(e')$ for all pairs $(e,e')$ of the form
$(e_{ij},e_{i-1,j})$ or $(e_{ij},e_{i+1,j+1})$; we denote the set of
these pairs by $\Pi$.

In what follows, when proving that one or another instance
$\Lscr$ of valid inequalities is facet-determining, we
will try to construct $2n+2m-2$ GT-patterns for $G$ such that their
borders are linear independent and each of them turns $\Lscr$ into
equality. First of all we show that the cone $\Bstar$ has co-dimension
one (recall that $\Bstar$ is contained in the hyperplane $\Hscr$).

\medskip
\noindent {\bf Claim.}
{\it $\Bstar_{n,m}$ contains $2n+2m-1$ linearly independent vectors
(assuming $n\ge 1$ and $m\ge 0$).}

\medskip
\noindent {\bf Proof.}
Consider the $2n+2m-1$ functions $h_1,\ldots,h_{n+m},\bar h_1,
\ldots,\bar h_m,q_1,\ldots,q_{n-1}$ on $\Estar$ defined by
  $$
 h_j:=\chi^{e_j}; \qquad \bar h_j:=\chi^{\bar e_j}; \qquad
    q_i:=\chi^{e_{i,i+m}}.
  $$
Observe that the border of each $h_j$ ($\bar h_j$) takes value 1 on
the edge $e_j$ (resp. $\bar e_j$), and 0 on the other edges in
$L\cup\bar L\cup R$, while the border of $q_i$ takes value 1 on
$r_{i-1}$, $-1$ on $r_i$, and 0 otherwise. Therefore, these $2n+2m-1$
borders are linearly independent. However, the above functions, except
for $h_1$, are not GT-patterns (e.g., $h_2(e_{n-1,1})=0<1=h_2(e_2)$).
By this reason, add $ch$ to each of these functions, where $c$ is a
large positive number and $h$ is the function on $\Estar$ defined by
   $$
  h(e_{ij}):=i-2j.
   $$
One can see that $h$ is a GT-pattern; moreover, $h(e)-h(e')=1$ for
each pair $(e,e')\in\Pi$. This implies that the updated functions are
already GT-patterns. Also their borders remain linearly independent
(as $c$ is large). So we have $2n+2m-1$ linearly independent vectors
in $\Bstar$, as required. \qed

\medskip
Next consider $I\subseteq\{1,\ldots,n\}$ and $J\subseteq\{1,\ldots,
m\}$ with $0<|I|+|J|<n+m$. (When $|I|=n$ and $|J|=m$, \refeq{redhorn}
becomes $|\lambda|-|\bar\lambda|-|\nu|\ge 0$, which is not essential.)
We examine the three cases of $I$ indicated in the theorem.

\medskip
{\em Case 1}. Let $0<|I|<n$. We have to show that~\refeq{redhorn}
determines a facet for any $J$. Our construction of $2n+2m-2$
linearly independent vectors that belong to $\Bstar$ and
attain equality in~\refeq{redhorn} is based on a certain partition
of $\Estar$ into subsets $A_1,\ldots,A_{n+m}$ (depending on $I$). Each
$A_j$ is defined to be a minimal set satisfying the following
conditions:
 \begin{myitem1}
 \item[(i)] $A_j$ contains $e_j$;
 \item[(ii)] if $A_j$ contains $e_{ip}$ with $0\ne i\not\in I$ and
$p<i+m$, then $e_{i-1,p}\in A_j$;
 \item[(iii)] if $A_j$ contains $e_{ip}$ with $i\in I$ and $p>1$,
then $e_{i-1,p-1}\in A_j$.
  \label{eq:p_Aj}
  \end{myitem1}

These sets are defined uniquely and do give a partition $P$ of $\Estar$.
(For example: if $n=3$, $m=2$ and $I=\{2\}$, then
$A_1=\{e_1,e_{21}\}$, $A_2=\{e_2,e_{22},e_{11},\bar e_1\}$, $A_3=\{e_3,
e_{23},e_{12},\bar e_2\}$, $A_4=\{e_4,e_{24},e_{13}\}$,
$A_5=\{e_5\}$.) Also for $j=1,\ldots,m$, both edges
$\bar e_j,e_{j+|I|}$ are contained in the same set in $P$, namely, in
$A_{j+|I|}$.

Define the function $h$ on $\Estar$ by
  \begin{equation} \label{eq:f_h}
h(e):=n+m-j \qquad \mbox{for $e\in A_j$, $j=1,\ldots,n+m$.}
  \end{equation}
The following properties of $h$ will be important for us:
  \begin{myitem1}
 \item[(i)] $h$ is a GT-pattern, and $h(e)>h(e')$ holds for each pair
$(e,e')\in\Pi$ such that $e,e'$ belong to different sets in $P$;
 \item[(ii)] \refeq{redhorn} holds with equality for the border
$(\lambda,\bar\lambda,\nu)$ of $h$.
  \label{eq:pr_h}
  \end{myitem1}
Indeed, (i) follows from~\refeq{f_h} because if $e_{ip}\in A_j$ then
each of $e_{i-1,p}$ and $e_{i+1,p+1}$ (if any) belongs to $A_j$ or
$A_{j+1}$, by~\refeq{p_Aj} in the construction of $P$. To see (ii), let
$\ell_j:=n+m-j$ for $j=1,\ldots,n+m$. Then $\lambda[1,|I|]=\ell_1+
\ldots+\ell_{|I|}$ and $\lambda(J+|I|)=\sum(\ell_j: j\in J+|I|)=
\bar\lambda(J)$. Observe that for each $i\in I$ and $p=2,\ldots,i+m$,
both edges $e_{ip},e_{i-1,p-1}$ belong to the same set in $P$, and
the edge $e_{i,1}$ belongs to the set $A_{d(i)}$, where $d(i)$ is the
number of elements of $I$ greater than or equal to $i$. This implies
  $$
 \nu_i=h(e_{i,1})+\sum_{p=2}^{i+m}(h(e_{ip})-h(e_{i-1,p-1})
          =\ell_{d(i)},
  $$
whence $\nu(I)=\sum(\ell_{d(i)}:i\in I)=\ell_1+\ldots+\ell_{|I|}$.
Therefore, $\lambda[1,|I|]=\nu(I)$ and $\lambda(J+|I|)=
\bar\lambda(J)$, yielding\refeq{pr_h}(ii).

Next we construct $2n+2m-2$ functions on $\Estar$, not necessarily
GT-patterns, such that their borders are linearly independent and
attain equality in~\refeq{redhorn}. To this aim, we form two
auxilliary \trap-grids $G'=(V',E')$ and $G''=(V'',E'')$, the former
having size $(n':=|I|,m':=|J|)$ and the latter having size
$(n'':=n-n',m'':=m-m')$. For convenience, we use notation with primes
(double primes) for edges and their sets in $G'$ (resp. $G''$). Since
$0<|I|<n$, we have $n',n''>0$ (while $m',m''\ge 0$). By the Claim
applied to $G'$, there exist GT-patterns $a_1,\ldots,a_{2n'+2m'-1}:
\Epstar\to \Rset$ (concerning $G'$) whose borders are linearly
independent. Similarly, there exist GT-patterns
$b_1,\ldots,b_{2n''+2m''-1}:\Eppstar\to\Rset$ (concerning $G''$) whose
borders are linearly independent.

These patterns are transformed (``lifted'') into functions on $\Estar$
by use of special maps $\omega_1:E_1\to\Epstar$ and $\omega_2:E_2\to
\Eppstar$, where $E_1$ ($E_2$) is the union of the sets $A_j$ for
$j\in J_1:=\{1,\ldots,|I|\}\cup(J+|I|)$ (resp. for
$j\in J_2:=\{1,\ldots,n+m\}\setminus J_1$). The map $\omega_1$ is
defined so as to satisfy the following condition:
  \begin{myitem}
for $e_{ij},e_{pq}\in E_1$, $e'_{i'j'}:=\omega_1(e_{ij})$ and
$e'_{p'q'}=\omega_1(e_{pq})$,
  \vspace{-5pt}
  \begin{itemize}
 \item[(i)] if $i=p$ and $j<q$, then $i'=p'$ and $j'<q'$;
 \item[(ii)] whenever $e_{ij},e_{pq}$ belong to the same set $A_d$
($d\in J_1$): if $p=i-1$ and $q=j$ (and therefore, $i\not\in I$,
by~\refeq{p_Aj}) then $i'j'=p'q'$; and if $p=i-1$ and $q=j-1$ (and
therefore, $i\in I$) then $p'=i'-1$ and $q'=j'-1$.
  \end{itemize}
  \label{eq:omega1}
  \end{myitem}

In its turn, $\omega_2$ is defined so as to satisfy:
  \begin{myitem}
for $e_{ij},e_{pq}\in E_2$, $e''_{i''j''}:=\omega_2(e_{ij})$ and
$e''_{p''q''}=\omega_2(e_{pq})$,
  \vspace{-5pt}
  \begin{itemize}
 \item[(i)] if $i=p$ and $j<q$, then $i''=p''$ and $j''<q''$;
 \item[(ii)] whenever $e_{ij},e_{pq}$ belong to the same set $A_d$
($d\in J_2$): if $p=i-1$ and $q=j$, then $p''=i''-1$ and $q''=j''$;
and if $p=i-1$ and $q=j-1$, then $i''j''=p''q''$.
  \end{itemize}
  \label{eq:omega2}
  \end{myitem}

One can check that both $\omega_1,\omega_2$ are well-defined and
unique. Also $\omega_1$ establishes one-to-one correspondence between
the sets $L_1:=\{e_j: j\in J_1\}$ and
$L':=\{e'_1,\ldots,e'_{n'+m'}\}$, as well as between the sets
$\bar L_1:=\{\bar e_j: j\in J\}$ and $\bar L':=\{\bar e'_1,\ldots,
\bar e'_{m'}\}$. Similarly, $\omega_2$ establishes one-to-one
correspondence between $L_2:=\{e_j: j\in J_2\}$ and
$L'':=\{e''_1,\ldots,e''_{n''+m''}\}$, and between
$\bar L_2:=\bar L\setminus \bar L_1$ and
$\bar L'':=\{\bar e''_1,\ldots,\bar e''_{m''}\}$.

Using $\omega_1$ and $\omega_2$, the above-mentioned functions
$a_s$ and $b_t$ are transformed into functions on $\Estar$ in a
natural way, as follows. For $s=1,\ldots,n'+m'-1$, define
  \begin{equation} \label{eq:f_gs}
 g_s(e):=\left\{
    \begin{aligned}
   a_s(\omega_1(e)) & \qquad \mbox{for $e\in E_1$}, \\
   0                & \qquad \mbox{for $e\in E_2$},
    \end{aligned}
          \right.
  \end{equation}
and for $t=1,\ldots,2n''+m''-1$, define
  \begin{equation} \label{eq:f_ht}
 h_t(e):=\left\{
    \begin{aligned}
   b_t(\omega_1(e)) & \qquad \mbox{for $e\in E_2$}, \\
   0                & \qquad \mbox{for $e\in E_1$}.
    \end{aligned}
          \right.
  \end{equation}

Let $\sigma_s$ ($\zeta_t$) be the border of $g_s$ (resp. $h_t$). The
values of $\sigma_s$ on the edges in $L_1\cup\bar L_1$ coincide with
the values of the border $\alpha_s$ of $a_s$ on the corresponding
edges in $L'\cup \bar L'$, and $\sigma_s$ takes zero values on the
remaining edges in $L\cup\bar L$.
Also~\refeq{omega1} and~\refeq{f_gs} show that for
$i=0,\ldots,n$, the sums $g_s(e_{i,1})+\ldots+g_s(e_{i,i+m})$ and
$a_s(e'_{i',1})+\ldots+a_s(e'_{i',i'+m'})$ are equal, where
$i'$ is the number of elements of $I$ smaller than or equal to $i$.
This implies that the value of $\sigma_s$ on the edge $r_i$ of the
right boundary of $G$ is zero if $i\not\in I$, and equals the value of
$\alpha_s$ on the edge $r'_{i'}$ of the right boundary of $G'$ if
$i\in I$. Similar correspondences take place for $\zeta_t$ and the
border of $b_t$ (regarding the sets $L_2,\bar L_2,L'',\bar L''$ and
replacing $i'$ by $i'':=|\{0,\ldots,i\}\setminus I|$).

So we can conclude that these $2n+2m-2$ vectors $\sigma_s$ and
$\zeta_t$ are linear independent. Also the equality $|\lambda'|-
|\bar\lambda'|-|\nu|=0$ for each $\alpha_s=(\lambda',\bar\lambda',
\nu')$ implies $\lambda(J_1)-\bar\lambda(J)-\nu(I)=0$ for $\sigma_s=
(\lambda,\bar\lambda,\nu)$, while the latter equality holds
automatically for each $\zeta_t=(\lambda,\bar\lambda,\nu)$ since
$\zeta_t$ is zero within $L_1$, $\bar L_1$ and $\{r_i: i\in I\}$.
Thus, each of the obtained borders turns~\refeq{redhorn} into
equality.

Since the functions $g_s$ and $h_t$ need not be GT-patterns, we add
$ch$ to each of them, where $h$ is defined in~\refeq{f_h} and $c$ is a
large positive number. Then the new functions are GT-patterns and
their borders are as required. (The former relies on
property~\refeq{pr_h} and the relationships between $g_s,\zeta_t$ and
$a_s,b_t$: observe that for $(e,\tilde e)\in\Pi$, if both $e,\tilde e$
are in the same set of the partition $P$, then $g_s(e)\ge
g_s(\tilde e)$, $h_t(e)\ge h_t(\tilde e)$ and $h(e)=h(\tilde e)$, while
if they are in different sets, then $h(e)>h(\tilde e)$.)
This completes the proof for case (i) in the theorem.

\medskip
{\it Case 2}. Let $I=\emptyset$. Then~\refeq{redhorn} becomes
$\lambda(J)-\bar\lambda(J)\ge 0$, which is the sum of valid
inequalities $\lambda_j-\bar\lambda_j\ge 0$ over $j\in J$. So we have
to examine only the latter inequalities for $j=1,\ldots,m$ (each
concerning the case $|J|=1$). To show that each
$\lambda_j-\bar\lambda_j\ge 0$ is facet-determining, we argue as in
Case 1 and apply the reduction to the grids $G',G''$ for $J=\{j\}$.
The grid $G''$ has size $(n,m-1)$ and generates $2n+2(m-1)-1=2n+2m-3$
linearly independent vectors in the intersection of $\Bstar$ and the
hyperplane $\lambda_j-\bar\lambda_j=0$. Moreover, these vectors
$(\lambda,\bar\lambda,\nu)$ satisfy $\lambda_j=\bar\lambda_j=0$. One
more vector comes up from the degenerate grid $G'$ (having size
(0,1)); it has entries $\lambda_j=\bar\lambda_j=1$, and 0 otherwise.
(This is the border of the GT-pattern taking value 1 on the edges
$e_{0,j},e_{1,j},\ldots,e_{nj}$, and 0 on the remaining edges in
$\Estar$.) This yields (ii) in the theorem.

\medskip
{\it Case 3}. Let $|I|=n$. This case is symmetric to Case 2.
Inequality~\refeq{redhorn} becomes $\lambda[1,n]+\lambda(J+n)
-\bar\lambda(J)-|\nu|\ge 0$. This is equivalent to
$-\lambda(\bar J+n)+\bar\lambda(\bar J)\ge 0$, where $\bar J:=
\{1,\ldots,m\}\setminus J$ (assuming $|\lambda|-|\bar\lambda|
-|\nu|=0$). The latter is the sum of valid inequalities
$-\lambda_{j+n}+\bar\lambda_j\ge 0$ over $j\in\bar J$ (each being
equivalent to~\refeq{redhorn} with an $(m-1)$-element set as $J$).
To show that for each $j=1,\ldots,m$, the inequality
$-\lambda_{j+n}+\bar\lambda_j\ge 0$ determines a facet of $\Bstar$,
we apply the reduction to the grids $G',G''$ for $J=\{1,\ldots,m\}
\setminus \{j\}$. Then $2n+2m-3$ vectors come up from $G'$, and one
vector from $G''$. (Another method: use the symmetry of $\Bstar$
defined by $(\lambda,\bar\lambda,\nu)\to -(\lambda_{n+m},\ldots,
\lambda_1,\bar\lambda_m,\ldots,\bar\lambda_1,\nu_1,\ldots,\nu_n)$,
which reduces case (iii) to case (ii) in the theorem.)

\medskip
It remains to examine the chamber inequalities. When $n=1$, each
inequality $\lambda_j\ge\lambda_{j+1}$ is not essential, as it is the
sum of $\lambda_j\ge\bar\lambda_j$ and $\bar\lambda_j\ge\lambda_{j+1}$
(which are the instances of~\refeq{redhorn} with $J:=\{j\}$ and with
$J:=\{1,\ldots,m\}-\{j\}$), and similarly for the inequalities
$\bar\lambda_{j'}\ge\bar\lambda_{j'+1}$. When $n=2$ and $m=0$, the
unique chamber inequality $\lambda_1\ge\lambda_2$ is not essential as
well, as it follows from the inequality
$\lambda_1-\nu_1\ge 0$ (i.e.,~\refeq{redhorn} for $I:=\{1\}$), the
inequality $\lambda_1-\nu_2\ge 0$ (i.e.,~\refeq{redhorn} for
$I:=\{2\}$), and the equality $\lambda_1+\lambda_2-\nu_1-\nu_2=0$.

Now consider the case $n\ge 2$ and $n+m\ge 3$. We show that for $j=
1,\ldots,n+m-1$, $\lambda_j\ge\lambda_{j+1}$ determines a facet of
$\Bstar$, as follows. Take the following vectors:

the border of $\chi^{e_j}$ for $j=1,\ldots,j-1,j+2,\ldots, n+m$;

the border of $\chi^{\bar e_j}$ for $j=1,\ldots,m$;

the border of $\chi^{e_{ip}}$ for $i=1,\ldots,n-1$ and one $p$ such
that $1\le p\le i+m$ and $(i,p)\ne (n-1,j)$;

the border $\xi$ of the all unit function on $\Estar$.

\smallskip
\noindent (An edge $e_{n-1,p}$ different from $e_{n-1,j}$ exists since
$n+m\ge 3$.) These $2n+2m-2$ vectors are linear independent and
satisfy the equality $\lambda_j=\lambda_{j+1}$. We modify each of them
by adding $c\sigma$, where $c$ is a large positive number and $\sigma$
is the border of a GT-pattern $h$ on $\Estar$
such that $h(e_j)=h(e_{j+1})=h(e_{n-1,j})$ and $h(e)>h(e')$ for all
pairs $(e,e')\in\Pi$ except for $(e_j,e_{n-1,j})$ and
$(e_{n-1,j},e_{j+1})$. (To construct such a pattern is easy.) The
new vectors become the borders of GT-patterns and are as required. The
fact that each inequality $\bar\lambda_j\ge \bar\lambda_{j+1}$ (when
$n,m\ge 2$) determines a facet is proved in a similar way.

We leave it for the reader to verify that all facets of $\Bstar$
appeared throughout the above proof are indeed different.
(In other words, for the facet-determining inequalities claimed in
the theorem, their incidence $(0,\pm 1)$-vectors have the property
that no distinct vectors $\xi,\xi'$ satisfy $\xi=\alpha\xi'
+\beta\theta$, where $\alpha\in \Rset_+$, $\beta\in\Rset$, and
$\theta$ is the incidence vector of the equality
$|\lambda|-|\bar\lambda|-|\nu|=0$.)

This completes the proof of Theorem~\ref{tm:t3}.

\section{Proof of Theorem 2} \label{sec:proof2}

First of all we observe that the generic case of convex
configuration in this theorem is reduced to the case of
\trap-configuration. Indeed, given $\lambda,\bar\lambda,\mu$ for
${V}$ as in~\refeq{conv_arr}, there exists a (sufficiently large)
positive integer $c$ such that each face of $\SC(\lamblam,\mu)$
contains a face of the polyhedron $\Pscr$ formed by the arrays
$X\in\SC(\lamblam,\mu)$ with $|\partial x_{ij}|\le c$
for all entries $\partial x_{ij}$ of $\partial X$. Let $m:=b_0$
and extend $\lambda$ to $(n+m)$-tuple $\lambda'$ by setting
$\lambda'_1:=\ldots :=\lambda_{a_n}:=c$, $\lambda'_{b(n)+1}:=\ldots
:=\lambda'_{n+m}:=-c$ and $\lambda'_j:=\lambda_j$ for $j=a_n+1,
\ldots,b_n$. Accordingly, set $\mu'_i:=\mu_i$ if $a_i=0$, and
$\mu'_i:=\mu_i-c$ if $a_i>0$. Then the restriction map $X'\to
{X'}\rest{{V}}$ gives a bijection between the \trap-arrays $X'$ with
$\lambda^{X'}=\lambda'$, $\bar\lambda^{X'}=\bar\lambda$, $\mu^{X'}=
\mu'$ and the arrays in $\Pscr$ (cf. explanations in the
Introduction). This implies that $\SC(\lamblam,\mu)$ is integral if
$\SC(\lambda'\setminus\bar\lambda,\mu')$ is such.

In the rest of the proof we deal with \trap-configuration of size
$(n,m)$. As before, we may assume $\mu=0^n$. Also one may assume
that $\lambda$ is nonnegative (cf. reasonings in the previous
section). For brevity we denote the polytope $\SC(\lamblam,0^n)$ by
$\SC(\lamblam)$. Theorem~2 will be proved by
constructing a bijection between the vertices of
$\SC(\lamblam)$ and certain forests in the grid $G$
(defined in Remark~1 in the Introduction). Establishing
this correspondence, we admit $\lambda$ and $\bar\lambda$ to be
real-valued.

The node set $V$ of $G$ is naturally partitioned into subsets
({\em horizontal layers}) $L_i=\{(i,0),$ $\ldots,(i,i+m)\}$,
$i=0,\ldots,n$. Extract the edges connecting neighbouring layers and
orient them from the top to the bottom. Formally: let $A$ be the
set of pairs $e^0_{ij}:=((i,j),(i+1,j))$ and
$e^1_{ij}:=((i,j),(i+1,j+1))$ of nodes of $G$, for
$i=0,\ldots,n-1$, $j=0,\ldots,i+m$. Then
$H_{n,m}:=H:=(V,A)$ is an acyclic digraph in which any maximal
(directed) path begins at a node of the ``topmost'' layer $L_0$ and
ends at a node of the ``bottommost'' layer $L_n$.

We say that a function $g:A:\to\Rset_+$ is a
$(\lambda,\bar\lambda)$-{\em admissible flow} in $H$ if
\begin{equation}\label{eq:div}
   \diver_g(i,j) := \left\{
    \begin{array}{rl}
  0, &\qquad i=1,\ldots,n-1,\;\; j=0,\ldots,i+m,\\
  \lambda_j-\lambda_{j+1}, & \qquad i=n,\;\; j=0,\ldots,n+m,\\
  \bar\lambda_{j+1}-\bar\lambda_j, & \qquad i=0,\;\;j=0,\ldots,m.
  \end{array}
        \right.
   \end{equation}
Here $\diver_g(v)$ ($v\in V$) stands for the value
$\sum_{e=(u,v)\in A}g(e)-\sum_{e=(v,u)\in A}g(e)$, and we formally
extend $\lambda$ and $\bar\lambda$ by setting $\lambda_0:=
\bar\lambda_0:=\lambda_1$ and $\lambda_{n+m+1}:=\bar\lambda_{m+1}:=0$.
In particular, $g(e^0_{n-1,0})=0$ and $g(e^1_{n-1,n+m-1})=\lambda_{n+m}$.
  \Xcomment{
Note also that $g$ is determined by its values on the subset
$A':=A\setminus\{e^0_{00},e^0_{10},\ldots,e^0_{n-1,0},e^1_{00},
e^1_{1,m+1},\ldots,e^1_{n-1,n+m-1}\}$.
  }
The set $\Fscr(\lamblam)$ of
$(\lambda,\bar\lambda)$-admissible flows forms a polytope in
$\Rset^{|A|}$.

\bigskip
\noindent
{\bf Claim.} {\sl For any $X\in\SC(\lamblam)$ there
exists a $(\lambda,\bar\lambda)$-admissible flow $g=\gamma(X)$
satisfying
  \begin{equation}  \label{eq:ge}
    \begin{array}{rl}
  g(e^0_{ij})= &\partial x_{i,j}-\partial x_{i+1,j+1}, \\
  g(e^1_{ij})= &\partial x_{i+1,j+1}-\partial x_{i,j+1},
    \end{array}
\qquad i=0,\ldots,n-1,\;\; j=0,\ldots,i+m,
  \end{equation}
letting $\partial x_{i0}:=\lambda_1$ and $\partial x_{i,i+m+1}:=0$.
Moreover, $\gamma$ is a bijective mapping of
$\SC(\lamblam)$ to
$\Fscr(\lamblam)$. }

\begin{figure}[htb]
 \begin{center}
  \unitlength=1mm
  \begin{picture}(45,20)
\put(0,0){\circle*{1.5}}
\put(8,0){\circle*{1.5}}
\put(16,0){\circle*{1.5}}
\put(24,0){\circle*{1.5}}
\put(32,0){\circle*{1.5}}
\put(40,0){\circle*{1.5}}
\put(4,6){\circle*{1.5}}
\put(12,6){\circle*{1.5}}
\put(20,6){\circle*{1.5}}
\put(28,6){\circle*{1.5}}
\put(36,6){\circle*{1.5}}
\put(8,12){\circle*{1.5}}
\put(16,12){\circle*{1.5}}
\put(24,12){\circle*{1.5}}
\put(32,12){\circle*{1.5}}
\put(12,18){\circle*{1.5}}
\put(20,18){\circle*{1.5}}
\put(28,18){\circle*{1.5}}
\put(8,0){\line(-2,3){4}}
\put(24,0){\line(-2,3){4}}
\put(40,0){\line(-2,3){4}}
\put(20,6){\line(-2,3){4}}
\put(36,6){\line(-2,3){4}}
\put(24,12){\line(-2,3){4}}
\put(31.5,12){\line(-2,3){4}}
\put(32.5,12){\line(-2,3){4}}
\put(8,0){\line(2,3){4}}
\put(16,0){\line(2,3){4}}
\put(24,0){\line(2,3){4}}
\put(4,6){\line(2,3){4}}
\put(12,6){\line(2,3){4}}
\put(20,6){\line(2,3){4}}
\put(28,6){\line(2,3){4}}
\put(8,12){\line(2,3){4}}
\put(15.5,12){\line(2,3){4}}
\put(16.5,12){\line(2,3){4}}
  \end{picture}
 \end{center}
 \caption{the flow corresponding to $\partial X$ in
Fig.~\ref{fig:der}b.}
  \label{fig:flow}
  \end{figure}
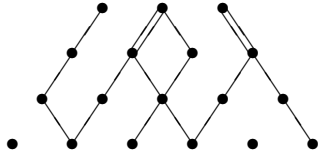

\medskip
(Figure~\ref{fig:flow} illustrates the flow determined by the array
$X$ with $\partial X$ as in Fig.~\ref{fig:der}b; here the flow is
integer and its value on an edge is indicated by the number of lines
connecting the ends of this edge.)

\medskip
\noindent
{\bf Proof.} Let $X\in\SC(\lamblam)$ and let $g$
be defined by~\refeq{ge}. Then for each node $v=(n,j)$ with
$j=0,\ldots,n+m$,
   $$
  \diver_g(v)=g(e^1_{n-1,j-1})+g(e^0_{n-1,j})=
\left(\partial x_{nj}-\partial x_{n-1,j}\right)+
\left(\partial x_{n-1,j}-\partial x_{n,j+1}\right)=
          \lambda_j-\lambda_{j+1},
  $$
letting $g(e):=0$ if the edge $e$ is void (e.g., for
$e=e^1_{n-1,-1}$).
Similarly, $\diver_g(v)=\bar\lambda_{j+1}-\bar\lambda_j$ for each
node $v=(0,j)$, $j=0,\ldots,m$. And for each node $v=(i,j)$
with $1\le i\le n-1$ and $0\le j\le i+m$, one has
  \begin{multline*}
\diver_g(v)=g(e^1_{i-1,j-1})+g(e^0_{i-1,j})-g(e^0_{ij})-g(e^1_{ij}) \\
  =\left(\partial x_{ij}-\partial x_{i-1,j}\right)
    +\left(\partial x_{i-1,j}-\partial x_{i,j+1}\right)
      -\left(\partial x_{ij}-\partial x_{i+1,j+1}\right)
        -\left(\partial x_{i+1,j+1}-\partial x_{i,j+1}\right) =0.
  \end{multline*}
Also the function $g$ is nonnegative, as is seen by
comparing~\refeq{ge} and~\refeq{concX}. Thus, $g$ is a
$(\lambda,\bar\lambda)$-admissible flow.
   \Xcomment{
(after the corresponding
extension to the edges in $A\setminus A'$).
   }

Conversely, let $g$ be a $(\lambda,\bar\lambda)$-admissible flow in
$H$. Assign numbers $\partial x_{ij}$ recursively by the following
rule:
   $$
    \begin{array}{rll}
  \partial x_{nj}:= & \lambda_j, & j=1,\ldots,n+m,\\
  \partial x_{ij}:= & \partial x_{i+1,j}-g(e^1_{i,j-1}),\qquad &
                   i=1,\ldots,n-1,\;\; j=1,\ldots,i+m.
  \end{array}
   $$
This gives the \trap-array $X$ of size $(n,m)$.
Reversing the argument above, one can
check validity of~\refeq{ge}. This and the nonnegativity of $g$ imply
that $X$ is strip-concave and satisfies $\lambda^X=\lambda$ and
$\bar\lambda^X=\bar\lambda$. Then $X\in\SC(\lambda\setminus
\bar\lambda)$, and the claim follows. \qed

\bigskip
Thus, $\gamma$ is a linear operator (in view of~\refeq{ge}) and
$\gamma$ gives a one-to-one correspondence between the points in the
polytopes $\SC(\lamblam)$ and
$\Fscr(\lamblam)$. Therefore, $\gamma$ establishes
a one-to-one correspondence between the vertices of these polytopes.

Next we characterize the vertices of
$\Fscr(\lamblam)$. To this aim,
we distinguish, in the bottommost layer $L_n$, the set $L(\lambda)$
of nodes $(n,j)$ ($1\le j\le n+m$) such that
$\lambda_j>\lambda_{j+1}$, and in the topmost layer $L_0$, the subset
$L(\bar\lambda)$ of nodes $(0,j)$ ($0\le j\le m$) such that
$\bar\lambda_j>\bar\lambda_{j+1}$.
Given a flow $g\in\Fscr(\lamblam)$, let
$H(g)$ denote the subgraph of $H$ induced by the set of edges $e$
with $g(e)>0$. From~\refeq{div} it follows that $H(g)$ contains
$L(\lambda)$ and $L(\bar\lambda)$ and that each node of $H(g)$ lies
on a path from $L(\bar\lambda)$ to $L(\lambda)$. Suppose there are two
different paths $P,P'$ in $H(g)$ having the same beginning and the
same end. Choose $\eps>0$ not exceeding the minimal value of $g$ on
the paths $P$ and $P'$. Then the functions
$g':=g+\eps\chi^P-\eps\chi^{P'}$ and $g'':=g-\eps\chi^P+\eps\chi^{P'}$
are nonnegative and satisfy~\refeq{div}, where $\chi^Q\in \{0,1\}^A$
is the characteristic function of the edge set of a path $Q$. So $g$
is expressed as the half-sum of two different
$(\lambda,\bar\lambda)$-admissible flows $g',g''$, and therefore,
$g$ cannot be a vertex of $\Fscr(\lamblam)$.

On the other hand, let for any two nodes $y$ and $z$, $H(g)$ contain
at most one path from $y$ to $z$,
i.e., $H(g)$ is a (directed) forest with the set $L(\bar\lambda)$ of
zero indegree nodes ({\em roots}) and the set $L(\lambda)$ of zero
outdegree ones ({\em leaves}). Then $g$ is the only
$(\lambda,\bar\lambda)$-admissible flow taking zero values on all edges
outside $H(g)$, i.e., $g$ is determined by $H(g)$. Indeed, one can see
that for each edge $e=(u,v)$ of $H(g)$, $g(e)$ is equal to
  \begin{equation} \label{eq:compon}
\sum(\lambda_j-\lambda_{j+1}:\; (n,j)\in V(Q))-
\sum(\bar\lambda_j-\bar\lambda_{j+1}:\; (0,j)\in V(Q)),
   \end{equation}
where $Q$ is the connected component of $H(g)\setminus\{e\}$ that
contains the node $v$, denoting by $V(Q)$ the node set of $Q$. This
implies that $g$ is a vertex of $\Fscr(\lamblam)$.
Moreover, $g$ is integer if $\lambda,\bar\lambda$ are integer, and
Theorem~\ref{tm:t2} follows. \qed

\smallskip
Arguing as in the above proof, one can associate the vertices of
$\SC(\lamblam)$ with certain subgraphs
of $H$, as follows.

  \begin{corollary} \label{cor:vert}
In case of \trap-configuration of size $(n,m)$,
each vertex of
$\SC(\lamblam)$ one-to-one corresponds to a forest
$H'$ in $H_{n,m}$ having $L(\bar\lambda)$ as the set of
roots and $L(\lambda)$ as the set of leaves and satisfying the
following condition: for each component $Q$ of $H'$, the value
in~\refeq{compon} is zero, and for each edge $e=(u,v)$ and the
component $Q$ of $H'\setminus\{e\}$ containing $v$, the value
in~\refeq{compon} is positive. Therefore, in case $m=0$,
the vertices of $\SC(\lambda)$ one-to-one correspond to the rooted
trees in $H_{n}$ ($:=H_{n,0}$) with root (0,0) and set of leaves
$L(\lambda)$.
  \end{corollary}

\noindent
\underline{Remark~6}.
The flows introduced in the proof of Theorem~2 give an alternative
way to represent the Gelfand-Tsetlin patterns (or the strip-concave
arrays), and Corollary~\ref{cor:vert} suggests a way to compute or
estimate the number of vertices of the polytope $\SC(\lambda
\setminus\bar\lambda,\mu)$ in case of \trap-configuration
(or $\Delta$-configuration). One can check that the reasonings in the
proof of Theorem~\ref{tm:t2} and the corresponding corollary are
applicable to \parall-configuration as well (with $H_{n,m}$ arising
from the corresponding parallelogram-wise grid of size $(n,m)$).

\section{Concluding Remarks} \label{sec:concl}

\medskip
In this section we outline (in parts A--D) more applications of the
flow approach developed in the proof of Theorem~\ref{tm:t2}. Here,
unless explicitly said otherwise, we consider the case of
\trap-configuration of size $(n,m)$. (Note that the exposed properties
remain valid if we deal with \parall-configuration.)

\medskip
{\bf A.}
Let $\Pscr=\Pscr_{n,m}$ be the set of paths in the graph $H=H_{n,m}$
beginning at a node of the layer $L_0$ and ending at a node of
$L_n\setminus\{(n,0)\}$. Associate with a path $P\in\Pscr$ the
\trap-array $Y^P$ with the entries $y_{i1}=\ldots=y_{i,p(i)}=1$ and
$y_{i,p(i)+1}=\ldots=y_{i,i+m}=0$ for $i=0,\ldots,n$, where
$(i,p(i))$ is a node of $P$. Considering the case of triangular
arrays, Berenstein and Kirillov~\cite{BK} noticed that the set of
arrays $Y^P$ ($P\in\Pscr_{n,0}$) constitutes a minimal list of
generators of the cone of nonnegative Gelfand-Tsetlin patterns of size
$n-1$. A similar property takes place for \trap-patterns (or
\parall-patterns) and can be easily shown by use of flows. More
precisely, for a strip-concave \trap-array $X$ with $\partial X\ge 0$,
take the flow $g=\gamma(X)$ defined by~\refeq{ge}. Then
$g$ is represented as a nonnegative linear combination
$\alpha_1\chi^{P_1}+\ldots+\alpha_N\chi^{P_N}$, where $P_1,\ldots,P_N
\in\Pscr$. One can check that $\partial X=\alpha_1Y^{P_1}+\ldots+
\alpha_NY^{P_N}$, as required (the minimality of $\{Y^P: P\in\Pscr\}$
is obvious).

\medskip
{\bf B.}
One can establish some invariants for polytopes
$\SC(\lamblam,0^n,\nu)$ when the entries of $\nu$
are permuted. Consider an array $X\in
\SC(\lamblam,0^n,\nu)$ and the flow $g=\gamma(X)$
as in~\refeq{ge}. For $i=1,\ldots,n$, we have
$\sum_{j=1}^{i+m}\partial x_{ij}-\sum_{j=1}^{i+m-1}\partial x_{i-1,j}=
x_{i,i+m}-x_{i-1,i+m-1}=\nu_i$. Also $\partial x_{ij}-
\partial x_{i-1,j}=g(e^1_{i-1,j-1})$ for $j=1,\ldots,i+m$
(see Section~\ref{sec:proof2} for the definition of edges
$e^0_{i'j'}$ and $e^1_{i'j'}$; as before, $\partial x_{i-1,i+m}:=0$).
Comparing these relations, we conclude that
  \begin{equation} \label{eq:nui}
\nu_i=g(e^1_{i-1,0})+\ldots+g(e^1_{i-1,i+m-1})\qquad
           \mbox{for $i=1,\ldots,n$.}
  \end{equation}

\begin{figure}[htb]
 \begin{center}
  \unitlength=1mm
  \begin{picture}(45,20)
\put(0,0){\circle*{1.5}}
\put(8,0){\circle*{1.5}}
\put(16,0){\circle*{1.5}}
\put(24,0){\circle*{1.5}}
\put(32,0){\circle*{1.5}}
\put(40,0){\circle*{1.5}}
\put(4,6){\circle*{1.5}}
\put(12,6){\circle*{1.5}}
\put(20,6){\circle*{1.5}}
\put(28,6){\circle*{1.5}}
\put(36,6){\circle*{1.5}}
\put(8,12){\circle*{1.5}}
\put(16,12){\circle*{1.5}}
\put(24,12){\circle*{1.5}}
\put(32,12){\circle*{1.5}}
\put(12,18){\circle*{1.5}}
\put(20,18){\circle*{1.5}}
\put(28,18){\circle*{1.5}}
\put(16,0){\line(-2,3){4}}
\put(40,0){\line(-2,3){4}}
\put(12,6){\line(-2,3){4}}
\put(28,6){\line(-2,3){4}}
\put(36,6){\line(-2,3){4}}
\put(24,12){\line(-2,3){4}}
\put(31.5,12){\line(-2,3){4}}
\put(32.5,12){\line(-2,3){4}}
\put(7.5,0){\line(2,3){4}}
\put(8.5,0){\line(2,3){4}}
\put(23.5,0){\line(2,3){4}}
\put(24.5,0){\line(2,3){4}}
\put(11.5,6){\line(2,3){4}}
\put(12.5,6){\line(2,3){4}}
\put(28,6){\line(2,3){4}}
\put(8,12){\line(2,3){4}}
\put(15.5,12){\line(2,3){4}}
\put(16.5,12){\line(2,3){4}}
  \end{picture}
 \end{center}
 \caption{the flow obtained by the zigzag swapping operation applied
to the flow in Fig.~\ref{fig:flow} at layer 2.}
  \label{fig:zig}
  \end{figure}
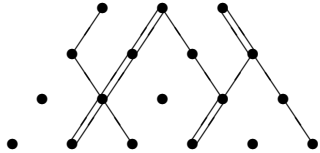

Choose $i\in\{1,\ldots,n-1\}$ and consider the subgraph $H^i$ of $H$
induced by the edges connecting the layers $L_{i-1},L_i$ or the
layers $L_i,L_{i+1}$. For $j=0,\ldots,i+m-1$, the nodes $(i-1,j)$
and $(i+1,j+1)$ are connected by two paths, namely, by path $Z_j$
with the edges $e^0_{i-1,j},e^1_{ij}$ and by path $Z'_j$ with the
edges $e^1_{i-1,j},e^0_{i,j+1}$. Let us call such a path $Z$ with
edges $e,e'$ a {\em zigzag} and define its capacity to be
$g(Z):=\min\{g(e),g(e')\}$. The {\em zigzag swapping operation}
modifies $g$ within $H^i$ by swapping the capacities simultaneously
for each pair $Z_j,Z'_j$. More precisely, for $j=0,\ldots,i+m-1$,
assign
       %
    $$
   g'(e) := \left\{
    \begin{array}{ll}
  g(e)-g(Z_j)+g(Z'_j)\qquad & \mbox{for each edge $e$ of $Z_j$},\\
  g(e)-g(Z'_j)+g(Z_j)\qquad & \mbox{for each edge $e$ of $Z'_j$},
  \end{array}
        \right.
   $$
  %
and $g'(e):=g(e)$ for the remaining edges of $H$. Obviously,
$g'$ is again a $(\lambda,\bar\lambda)$-admissible flow.
(For example, such an operation applied to the flow in
Fig.~\ref{fig:flow} results in the flow illustrated in
Fig.~\ref{fig:zig}.)

Let $\sigma_i(X)$ denote the array $\gamma^{-1}(g')$ and let $\nu'$
be the $n$-tuple of local differences on the ``right boundary'' of
$\sigma_i(X)$. Using~\refeq{nui}, one can check that the zigzag
swapping operation swaps $\nu_i$ and $\nu_{i+1}$, i.e.,
$\nu'_i=\nu_{i+1}$, $\nu'_{i+1}=\nu_i$ and $\nu'_p=\nu_p$ for
$p\ne i,i+1$. Moreover, applying the zigzag swapping operation (with
the same $i$) to $g'$ returns $g$.

Thus, for each $i$, $\sigma_i$ is a continuous bijective mapping of
$\SC(\lamblam,0^n,\nu)$ to $\SC(\lamblam,0^n,\nu')$
(and $\sigma_i^2$ is the identity on $\SC(\lamblam)$)%
   \footnote{
Note also that for integer points the zigzag swapping operation
produces Bender--Knuth's involution, cf.~\cite{BK}.
}.
Moreover, for $k\in\Nset$, if $g$ is $\frac1k$-integer, so is $g'$.
Therefore, $\sigma_i$ gives a bijection on the $\frac1k$-integer
points in these polytopes for any $k$. As a consequence (for $k=1$),
the following known property is
obtained: if $\lambda,\bar\lambda,\nu$ are integer and if $\nu'$ is an arbitrary
permutation of $\nu$, then Kostka coefficients $K(\lamblam,\nu)$ and
$K(\lamblam,\nu')$ are equal.

\medskip
{\bf C.}
Let $\lambda,\bar\lambda,\nu$ be rational-valued and let $\nu'$ be a
permutation of $\nu$. Let $V_0$ denote the set of boundary index pairs
in $V$ (or the boundary nodes in the grid $G$). The fact that each map
$\sigma_i$ is continuous and bijective implies that the polytopes
$\SC:=\SC(\lamblam,0^n,\nu)$ and $\SC':=\SC(\lamblam,0^n,\nu')$ have
the same dimension (which typically equals $|V\setminus V_0|$).
Consider the $|V\setminus V_0|$-dimensional affine subspaces
$S$ and $S'$ containing the polytopes $\SC$ and $\SC'$ , respectively,
which are obtained by imposing the corresponding equalities on the
values on $V_0$. Since $S$ and $S'$ are
parallel, there is $k'\in\Nset$ such that for any multiple $k$ of $k'$,
the lattice of $\frac1k$-integer points in $S'$ is obtained by a
parallel translation of a similar lattice in $S$. So the density of
$\frac1k$-integer points in $S$ and $S'$ (measured by the number of
such points in a unit ball with center at a point of the lattice) is
the same. Also the numbers of $\frac1k$-integer points in the
polytopes in question are equal. Thus, when $k$ tends to infinity, we
obtain equality for the corresponding volumes and can conclude with
the following.
  \begin{proposition}  \label{pr:permut}
Given (real-valued) $\lambda,\bar\lambda,\nu$,
let $\nu'$ be a permutation of $\nu$. Then the polytopes
$\SC(\lamblam,0^n,\nu)$ and $\SC(\lamblam,0^n,\nu')$ have the same
$|V\setminus V_0|$-dimensional volume.
  \end{proposition}

It should be noted that, although $\sigma_i$ (being a piece-wise
linear operator) brings integer points into integer ones, it need not
do so for polytope vertices, even for polytopes $\SC(\lamblam)$.
Indeed, in case $m=0$, take a rooted tree $T$ in $H_{n,0}$ (with root
(0,0) and the leaves in $L_n$) such that for some $i,j$, the subgraph
$T\cap H^i$ contains zigzags $Z_j$ and $Z'_{j+1}$. Then the zigzag
swapping operation (applied to a nowhere zero flow on $T$) transforms
the pair $Z_j,Z'_{j+1}$ into $Z'_j,Z_{j+1}$, so the resulting graph
$T'$ is not a tree, as it has two edges entering the node $(i,j+1)$.

\medskip
{\bf D.} The reduction applied in the proof of part ``if'' of
Theorem~\ref{tm:t1} in case $m>0$ can be described in terms of flows.
Moreover, the language of flows is convenient to develop a more
general sort of reduction and to demonstrate some additional
properties. To explain the idea, consider $X\in\SC(\lamblam,0^n,\nu)$
and $g$ as in~\refeq{ge}, assuming that $\lambda$ is nonnegative.
Let $P$ be a path in $H$ beginning at a node $(0,s)$ of the layer
$L_0$, ending at a node $(n,t)$ of the layer $L_n$ and such that the
minimum $\alpha$ of values of $g$ on the edges of $P$ is nonzero.
Choose $p\in\Zset$ and $\alpha'\in\Rset$ satisfying $0\le s+p\le m$ and
$0<\alpha'\le\alpha$ and change $g$ by moving the path $P$ with weight
$\alpha'$ at distance $|p|$, to the right of left depending on the
sign of $p$. Formally: define $P'$ to be the path containing the node
$(i,j+p)$ for each node $(i,j)$ of $P$ and transform $g$ into
$g':=g-\alpha'\chi^P+\alpha'\chi^{P'}$. This transformation does not
change the sum in~\refeq{nui}, and therefore, the
resulting array $X':=\gamma^{-1}(g')$ satisfies $\nu^{X'}=\nu$. When
$p>0$ ($p<0$), the row derivative $\partial X'$ is obtained from
$\partial X$ by increasing (resp. decreasing) by $\alpha'$ the entries
corresponding to the horizontal edges of the grid $G$ lying between
the paths $P$ and $P'$; the tuples $\lambda^X$ and $\bar\lambda^X$ are
changed accordingly.

Using such operations, one can transform $g$ more globally, still
preserving $\nu$: decompose $g$ into the sum of path flows
$\alpha_q\chi^{P_q}$ ($\alpha_q>0$), $q=1,\ldots,N$, and move each path
$P_q$ to the left so that the resulting $P'_q$ begin at the node (0,0).
This gives an array $X'$ with $\partial x'_{ij}=0$ for $i=0,\ldots,n$
and $j=i+1,\ldots,i+m$, i.e., in essense, $X'$ is equivalent to a
$\Delta$-array. One can deduce that the first $n$ entries of the tuple
$\lambda':=\lambda^{X'}$ are expressed as follows:
  \begin{equation}  \label{eq:transform}
\lambda'_k=\sum_{t=k}^{n+m} |[\lambda_t,\lambda_{t+1}]\cap
[\lambda_1,\bar\lambda_{t-k+1}]| \qquad \mbox{for $k=1,\ldots,n$},
  \end{equation}
denoting by $|[a,b]|$ the length $b-a$ of a segment $[a,b]$ and letting
$\bar\lambda_j:=0$ for $j>m$. Conversely, given $\lambda,\bar\lambda,
\nu$, define the $n$-tuple $\lambda'$ by~\refeq{transform} and
consider a $\Delta$-array $X'\in\SC(\lambda',0^n,\nu)$. Then one can
determine a special path decomposition for $\gamma(X')$ and move each
path at a due distance to the right so as to obtain a flow determining
a \trap-array in $\SC(\lamblam,0^n,\nu)$ (moreover, $\lambda'$ is
integer when $\lambda,\bar\lambda$ are such and one can maintain flow
and array intergality under the transformation). This gives a
constructive way to reduce the trapezoidal case to the triangular one.
The tuple $\lambda'$ is weakly decreasing and it just represents the
vertex generating vector for the permutohedron mentioned in the
Introduction.

\medskip
Next we explain the idea of deriving Theorem~\ref{tm:t1}
from results in~\cite{KT,KTW} (mentioned in Remark~2 in
Section~\ref{sec:intr}). We use the equivalence between \trap-arrays
of size $(n,m)$ and functions on the node set of the corresponding
grid $G=(V,E)$. Given tuples
$\lambda,\bar\lambda,\mu,\nu$, let us choose a positive integer $c$
and replace $\mu,\nu$ by $\mu',\nu'$ defined by $\mu'_i:=\mu_i-ic$ and
$\nu'_i:=\nu_i-ic$, $i=1,\ldots,n$. This turns the polytope
$\SC(\lamblam,\mu,\nu)$ into
$\SC(\lamblam,\mu',\nu')$ (each array $X$ in the
former polytope corresponds to $X'$ defined by $x'_{ij}:=x_{ij}-
\frac{i(i+1)}{2}c$); for brevity, we denote the latter polytope by
$\Cscr$. When $c$ is large enough, $\Cscr$ consists of fully concave
arrays, and we can apply results on the corresponding discrete concave
functions. The second part of Theorem~\ref{tm:t1} follows from a
result in~\cite{KT} (in fact, shown there for any convex grid) which
in our case reads: if $\lambda,\bar\lambda,\mu',\nu'$ are integer and
if $\Cscr\ne\emptyset$, then $\Cscr$ contains an integer point.

The first part of Theorem~\ref{tm:t1} follows from a combinatorial
characterization for the existence of a discrete concave function
under prescribed boundary data (we use its extension to an arbitrary
convex grid given in~\cite{Kar-03}). It uses a notion of {\em puzzle}
(originally introduced for $\Delta$-grids in~\cite{KTW}). This is a
subdivision $\Pi$ of the grid into a set of little triangles and
little rhombi (the union of two little triangles sharing an edge),
along with a 0,1-labeling of the edges of $G$ occurring in the
boundaries of these pieces, satisfying the following properties:

(i) for each little triangle $\tau$ in $\Pi$, the edges of $\tau$ are
all labeled either by 0 or by 1;

(ii) for each little rhombus $\rho$ in $\Pi$, a side edge of $\rho$ is
labeled 1 if clockwise of an obtuse angle, and 0 if clockwise of an
acute angle.

Then a necessary and sufficient condition on the non-emptiness of
$\Cscr$ (in \trap-case) is that each puzzle $\Pi$ satisfies the
inequality
  \begin{equation} \label{eq:puzzle}
\lambda(I)-\bar\lambda(J)+\mu'(K)-\nu'(L)\ge 0,
  \end{equation}
where $I,J,K,L$ are the sets of edges labeled 1 in the lower, upper,
left and right sides of $G$, respectively. To show the necessity is
rather easy, as follows. Let $\Cscr\ne\emptyset$ and let
$x\in\Cscr$ (considering $x$ as a function on $V$).
The discrete concavity of $x$ implies that for each little
rhombus $\rho$ with obtuse vertices
$u,u'$ and acute vertices $v,v'$, one has $q(x,\rho):=x(u)+x(u')
-x(v)-x(v')\ge 0$. When summing up these inequalities for all rhombi
in $\Pi$ and the equalities $(x(v)-x(u))+(x(w)-x(v))+(x(u)-x(w))=0$
for all little triangles labeled 1, with vertices $u,v,w$ in the
anticlockwise order, the terms $x(\cdot)$ for interior vertices
cancel out and we just obtain~\refeq{puzzle} with $I,J,K,L$ to be the
sets of edges labeled 1 on the corresponding sides.

When $c$ tends to $+\infty$, the value $q(x,\rho)$ does so as well
(uniformly for all $x\in\Cscr$) for each little rhombus $\rho$,
if any, whose smaller diagonal is parallel to the bottom side of $G$.
The grow of $q(x,\rho)$ must cause a similar behavior for the left
hand side in~\refeq{puzzle}. This implies that the puzzles containing
at least one of such rhombi $\rho$ can be excluded from the
consideration, as they become redundant in verification of the
non-emptiness of $\Cscr$. Now relation~\refeq{horn} in
Theorem~\ref{tm:t1} can be deduced from~\refeq{puzzle} when the
remaining puzzles $\Pi$ are considered.

\smallskip
In conclusion, it should be noted that, using the above reduction to
the fully concave case and an argument in~\cite{Bu} (where an
alternative proof of the integrality theorem from~\cite{KT} is given),
one can show the following sharper version of the last claim in
Theorem~\ref{tm:t1}.
  \begin{proposition}  \label{pr:downh}
For integer $\lambda,\bar\lambda,\mu,\nu$, the down hull $\Dscr$ of
$\SC(\lamblam,\mu,\nu)$ (i.e., the polyhedron $\SC(\lamblam,\mu,\nu)-
\Rset^{V}_+$) is integral.
   \end{proposition}

One can give a direct, relatively simple, proof of this proposition.
A sketch: Consider a vertex $X$ of $\Dscr$; then there is no array
$X'\ne X$ in $\Dscr$ with $X'\ge X$. Let ${V}_1,\ldots,{V}_N$ be
the minimal nonempty sets of index pairs such that for $q=1,\ldots,
N$ and for any $ij$ and $i'j'$ with $i'=i+1$, $j'\in\{j,j+1\}$ and
$\partial x_{ij}=\partial x_{i'j'}$, the set ${V}_q$ contains either
both or none of $ij$ and $i'j'$. Let $c_q:=\partial x_{ij}$ for
$ij\in{V}_q$. Each
${V}_q$ is associated with the corresponding subset of horizontal
edges in the grid $G$; let $R_q$ denote the union of
little triangles containing an edge in this subset.
Then the interior of each region $R_q$ is connected, and each maximal
horizontal line $\Lscr_i$ in $G$ (corresponding to the $i$-th row in
$\partial X$) intersects $R_q$ by a connected, possibly empty, set.
We say
that $R_q$ is an {\em intermediate} region if it has no edge in the
lower or upper boundary of $G$; let for definiteness $R_1,\ldots,
R_\ell$ be the intermediate regions. One shows that if the set
$V_q$ of nodes of $G$ occurring in the interior of an intermediate
region $R_q$ is nonempty, then one can increase the function $x$ by a
(small) positive constant within the set $V_q$ so as to preserve the
strip-concavity; the boundary tuples $\lambda^X,\bar\lambda^X,\mu^X,
\nu^X$ are preserved automatically. (This relies on the observation
that if, e.g., $\partial x_{ij}=\partial x_{i-1,j}$ and the vertex
$(i,j-1)$ is in $V_q$, then $(i-1,j-1)$ is in $V_q$ as well, in view of
$\partial x_{ij}=\partial x_{i,j-1}=\partial x_{i-1,j-1}$.)
Therefore, $V_q=\emptyset$ for
all $q=1,\ldots,\ell$; in other words, each horizontal line $\Lscr_i$
contains at most one edge within $R_q$.

Now associate with $R_q$ ($1\le q\le\ell$) a real variable $z_q$.
Let $A=(a_{iq})$ be the $(n-1)\times\ell$ matrix in which $a_{iq}$ is
the number of edges of the line $\Lscr_i$ occurring in $R_q$. Form
the linear system $Az=b$, where for $i=1,\ldots,n-1$, $b_i$ is equal to
$x_{i,i+m}-x_{i0}$ minus the sum of values $\partial x_{ij}$ over all
$ij$ concerning the edges of non-intermediate regions.
Then for the numbers $c_q$ as above,
the tuple $z:=(c_1,\ldots,c_\ell)$ is a solution to
this system. Note that each $b_i$ is an integer. (Indeed, each of the
above values $\partial x_{ij}$ is equal to some entry of $\lambda$ or
$\bar\lambda$, which is an integer; $x_{i0}$ and $x_{i,i+m}$ are
integers as well.) Also $A$ is a 0,1-matrix and the ones in each
column go in succession, i.e. $A$ is an {\em interval} matrix. So $A$
is totally unimodular (cf.~\cite[Section~19.4]{Schr}) and must have
full column rank (otherwize $Az=0$ has a nonzero solution and we can
represent $X$ as the half-sum of two other points in
$\SC(\lamblam,\mu,\nu)$). Then $c_1,\ldots,c_\ell$ are integers, as
required.

\medskip
{\bf Acknowledgement.} We thank the anonymous referee for correcting
some inaccuracies and suggesting improvements in the earlier versions
of this paper.

\end{document}